\documentclass[10pt]{amsart}
\usepackage{amsmath}
\usepackage{amsfonts}
\usepackage[usenames,dvipsnames,svgnames,table]{xcolor}
\usepackage{graphicx}
\usepackage[percent]{overpic}
\usepackage{amsthm}
\usepackage{amssymb}
\begingroup
    \makeatletter
    \@for\theoremstyle:=definition,remark,plain\do{%
        \expandafter\g@addto@macro\csname th@\theoremstyle\endcsname{%
            \addtolength\thm@preskip\parskip
            }%
        }
\endgroup

\usepackage[margin=10pt,font=small,labelfont=bf]{caption}

\usepackage[pdfauthor=Lyons\ and\ White,bookmarksnumbered,hyperfigures,pdftex,pagebackref,colorlinks=true,linkcolor=BrickRed,urlcolor=BrickRed,citecolor=OliveGreen,pdfstartview=FitH]{hyperref}

\usepackage{microtype}

\usepackage{mathtools}
\mathtoolsset{centercolon}

\numberwithin{equation}{section}
\numberwithin{figure}{section}

\usepackage{GWmacros}

\newcommand{\pa}{\mathcal{P}}
\renewcommand{\Pr}{\P}
\newcommand{\s}{x}
\newcommand{\e}{\varnothing}
\newcommand{\hm}{\nu_{\mathrm h}}
\newcommand{\hp}{\nu_+}
\newcommand{\sm}{\nu_{\mathrm s}}
\newcommand{\hmd}{\dual{\hm}}
\newcommand{\smd}{\dual{\sm}}
\newcommand{\bp}{o}
\newcommand{\dual}[1]{#1^\dagger}
\newcommand{\conc}{\oplus}
\newcommand\sfrac[2]{#1/#2}
\renewcommand{\defn}[1]{\textbf{\textit{#1}}}
\DeclareMathOperator{\esc}{Cr}
\newcommand\depth{\mathbf{D}}
\newcommand\states{\mathbf{\Gamma}}
\newcommand\rstates{\mathbf{\Gamma}_{{\bullet}}}
\newcommand\rstatesd{\mathbf{\Gamma}^\dagger_{{\bullet}}}
\newcommand\w{\Gamma_{\circ}}
\newcommand\wdd{\Gamma^\dagger_{\circ}}
\newcommand\Td{T^\dagger}
\newcommand\BS{\mathrm{BS}}

\newcommand\I[1]{\mathbf1_{#1}}
\newcommand\xor{\mathbin{\triangle}}

\newcommand\RL[1]{#1}
\let\P\relax
\DeclareMathOperator{\P}{\mathbf{P}\mathopen{}}
\DeclareMathOperator{\E}{\mathbf{E}\mathopen{}}

\def\Psub_#1{\P_{\! #1}}

\def\Pbig#1{\P\mkern-.5mu\bigl[#1\bigr]}
\def\Psubbig_#1#2{\Psub_{#1}\mkern-1.5mu\bigl[#2\bigr]}

\def\Psubbigg_#1#2{\Psub_{#1}\mkern-1.5mu\biggl[#2\biggr]}
\def\PBig#1{\P\mkern-.5mu\Bigl[#1\Bigr]}

\def\Ebig#1{\E\mkern-1.5mu\bigl[#1\bigr]}
\def\Esubbig_#1#2{\E_{#1}\mkern-1.5mu\bigl[#2\bigr]}
\def\EsubBig_#1#2{\E_{#1}\mkern-1.5mu\Bigl[#2\Bigr]}
\def\Esupbig^#1#2{\E^{#1}\mkern-1.5mu\bigl[#2\bigr]}
\def\EsupBig^#1#2{\E^{#1}\mkern-1.5mu\Bigl[#2\Bigr]}

\def\Esubbigg_#1#2{\E_{#1}\mkern-1.5mu\biggl[#2\biggr]}

\def\thmenv#1#2#3{\begin{#1} \label{#1:#2} #3 \end{#1}}
\def\richthmenv#1#2#3#4{\begin{#1}[#3] \label{#1:#2} #4 \end{#1}}

\def\procl#1.#2 #3\endprocl{%
       \ifx#1t\thmenv{Theorem}{#2}{#3}\fi
       \ifx#1l\thmenv{Lemma}{#2}{#3}\fi
       \ifx#1p\thmenv{Proposition}{#2}{#3}\fi
       \ifx#1c\thmenv{Corollary}{#2}{#3}\fi
       \ifx#1d\thmenv{Definition}{#2}{#3}\fi
       \ifx#1g\thmenv{Conjecture}{#2}{#3}\fi
       \ifx#1q\thmenv{Question}{#2}{#3}\fi
       \ifx#1r\thmenv{Remark}{#2}{#3}\fi
    }%

\def\rprocl#1.#2 #3 #4\endprocl{%
       \ifx#1t\richthmenv{Theorem}{#2}{#3}{#4}\fi
       \ifx#1l\richthmenv{Lemma}{#2}{#3}{#4}\fi
       \ifx#1p\richthmenv{Proposition}{#2}{#3}{#4}\fi
       \ifx#1c\richthmenv{Corollary}{#2}{#3}{#4}\fi
       \ifx#1d\richthmenv{Definition}{#2}{#3}{#4}\fi
       \ifx#1g\richthmenv{Conjecture}{#2}{#3}{#4}\fi
       \ifx#1q\richthmenv{Question}{#2}{#3}{#4}\fi
       \ifx#1r\richthmenv{Remark}{#2}{#3}{#4}\fi
    }%

\def\rref#1.#2/{%
      \ifx #1sSection~\ref{s.#2}\fi
      \ifx #1SSubsection~\ref{S.#2}\fi
      \ifx #1tTheorem~\ref{Theorem:#2}\fi  
      \ifx #1lLemma~\ref{Lemma:#2}\fi 
      \ifx #1cCorollary~\ref{Corollary:#2}\fi 
      \ifx #1pProposition~\ref{Proposition:#2}\fi 
      \ifx #1dDefinition~\ref{Definition:#2}\fi
      \ifx #1gConjecture~\ref{Conjecture:#2}\fi 
      \ifx #1qQuestion~\ref{Question:#2}\fi 
      \ifx #1rRemark~\ref{Remark:#2}\fi 
      \ifx #1aAppendix~\ref{a.#2}\fi 
      \ifx #1fFigure~\ref{f.#2}\fi
      \ifx #1e(\ref{e.#2})\fi
      \ifx #1b\cite{#2}\fi
      \ifx #1B\cite{#2}\fi
        }

\def\rlabel #1 #2{\begin{equation} \label{#1} #2 \end{equation}}

\def\rproof{\begin{proof}}

\def\eqaln#1{\begin{align*} #1 \end{align*}}

\def\rcases#1{\begin{cases} #1 \end{cases}}

\def\Qed{\end{proof}}

\def\bsection#1#2{\bigbreak\section{#1}\label{#2}}

\title[Dyadic Lattice Graphs]{A Stationary Planar Random Graph with\\ Singular Stationary Dual:
Dyadic Lattice Graphs}
\author{Russell Lyons and Graham White}
\address{Department of Mathematics, 831 E. 3rd St.,
Indiana University, Bloomington, IN 47405-7106} 
\thanks{%
The work of R.L.\ is partially supported by the National
Science Foundation under grant DMS-1612363.}

\email{\href{mailto:rdlyons@indiana.edu}{rdlyons@indiana.edu}}
\email{\href{mailto:grrwhite@iu.edu}{grrwhite@iu.edu}}

\begin{document}

\begin{abstract}
Dyadic lattice graphs and their duals are commonly used as discrete
approximations to the hyperbolic plane.
We use them to give examples of random rooted graphs that are stationary
for simple random walk, but whose duals have only a singular stationary measure.
This answers a question of Curien and shows behaviour different from the unimodular case.
The consequence is that planar duality does not combine well with stationary random graphs. 
We also study harmonic measure on dyadic lattice graphs and show its singularity.
\end{abstract}



\vglue -1cm
\maketitle


%
%

\section{Introduction}

Since the study \rref b.BLPS/ of group-invariant percolation, the use of unimodularity via the mass-transport principle has been an important tool in analysing percolation and other random subgraphs of Cayley graphs and more general transitive graphs; see, e.g., \cite[Chapters 8, 10, and 11]{LyonsPeres}. These ideas were extended by \rref b.AldousLyons/ to random rooted graphs. For planar graphs, a crucial additional tool is, of course, planar duality. In the deterministic case, \cite[Theorem 8.25]{LyonsPeres} shows that every planar quasi-transitive graph with one end is unimodular and admits a plane embedding whose plane dual is also quasi-transitive (hence, unimodular). This was extended by \cite[Example 9.6]{AldousLyons} to show that every unimodular random rooted plane graph \RL{satisfying a mild finiteness condition} admits a natural unimodular probability measure on the plane duals; in fact, the root of the dual can be chosen to be a face incident to the root of the primal graph. The recent paper \rref b.AHNR/ makes a systematic study of unimodular random planar graphs, synthesizing known results and introducing new ones, showing a dichotomy involving 17 equivalent properties. 

One significant implication of unimodularity is that when the measure is biased by the degree of the root, one obtains a stationary measure for simple random walk. This led Benjamini and Curien \rref b.BC:stat/ to study the general context of probability measures on rooted graphs that are stationary for simple random walk. One aim has been to elucidate which properties hold without the assumption of unimodularity. In particular, Nicolas Curien has asked (unpublished) whether stationary random graphs have stationary duals; our interpretation of this is the following question:

\begin{Question} \label{q.curien}
Given a probability measure $\mu$ on rooted plane graphs $(G, \bp)$ (that are locally finite and whose duals $\dual G$ are locally finite), let $\nu$ be the probability measure on rooted graphs $(\dual G, f)$ obtained from choosing a neighbouring face $f$ of $\bp$ uniformly at random. If $\mu$ is stationary (for simple random walk), then is $\nu$ mutually absolutely continuous to a stationary measure?
\end{Question}

In Question \ref{q.curien} and henceforth, we use `stationary' to mean `stationary with respect to simple random walk'. 

\begin{figure}[ht]
\begin{center}
\includegraphics[width=.6\textwidth]{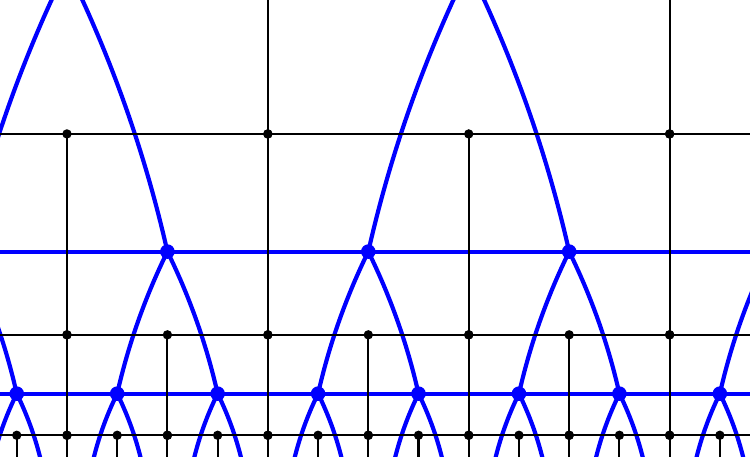}
\end{center}
\caption{A portion of the dyadic lattice graphs, with dual graph marked in blue.}
\label{fig:gammadual-intro}
\end{figure}

Here, we are actually interested in rooted plane graphs up to rooted isomorphisms induced by orientation-preserving homeomorphisms of the plane (though one could allow such a homeomorphism to change the orientation of the plane without affecting our results).
In this paper, we give a negative answer in the general stationary case. This means that planar duality does not combine well with stationary random graphs. Our counterexample uses dyadic lattice graphs; the primal graphs have vertices of only two different degrees and the dual graphs are regular. See Figure \ref{fig:gammadual-intro} for a representation of a portion of the primal and dual graphs.

On the other hand, if the primal and dual graphs are both regular, then the resulting graphs are uniquely determined by their degrees and codegrees, are transitive, and are unimodular: see page 197 and Theorem 8.25 of \cite{LyonsPeres}. In this sense, our counterexample is the best possible. See also the discussion of vertex and face degrees in Section \ref{sec:optimality}.
\RL{We do not have any examples, other than unimodular ones, of a stationary random plane graph with stationary dual; it seems possible there are none.}

The dyadic lattice graphs we study are often used as a combinatorial approximation of the hyperbolic
plane; see, e.g., the introductory work \cite[Section 14]{CFKP}.
Dyadic lattice graphs are also closely related to Whitney decompositions,
which have been used in studying diffusions since the work of Ba\~nuelos \rref b.Banuelos/.
Finally, these graphs are subgraphs of the usual Cayley graph of the
Baumslag--Solitar group $\BS(1,2)$---see Remark \ref{rem:bs}.

It might appear that there is only one dyadic lattice graph and only one dual. This is not true. While there is only one way to subdivide going downwards one level in the picture, there are two ways to agglomerate going upwards one level. Therefore, there are uncountably many such graphs. We will show, however, that there is a unique probability measure on rooted versions of such graphs that is stationary for simple random walk; the same holds for the duals. With appropriate notation, it is easy to put these measures on the same space; our main theorem is that these measures are mutually singular (\rref t.curien/). We also show that simple random walk tends downwards towards infinity and defines a harmonic measure. We will show that this measure is singular with respect to Lebesgue measure in a natural sense (Proposition \ref{prop:mutuallysingular}).

We will not define `unimodular' here because we will not use it again.
In Section \ref{sec:notation}, we give crucial notation \RL{for} the vertices in dyadic lattice graphs. This will also enable us to give useful notation to the entire rooted graph, which will identify a rooted dyadic lattice graph with a dyadic integer. That section also contains basic properties of dyadic lattice graphs that are invariant under automorphisms. In \rref s.basic/, we prove the fundamental existence and uniqueness properties of stationary and harmonic measures. In Section \ref{sec:reflections}, we show how various symmetries of dyadic lattice graphs lead to (sometimes surprising) comparisons and identities for random-walk probabilities. We then use these to prove the singularity results mentioned above. We do not have explicit formulas for either the stationary measure on the primal graphs nor the harmonic measure for the primal graphs. Thus, in Section \ref {sec:numerics}, we present some numerical approximations to these. Finally, Section \ref{sec:optimality} contains further discussion of the optimality of our example.

Sections \ref{sec:reflections} and \ref{sec:numerics} suggest several open questions. In particular, we do not know how to determine the stationary measure of even the simplest sets, or how to explain the patterns in the harmonic measure illustrated in Section \ref{sec:numerics}.

\medbreak
\RL{\textbf{Acknowledgements.} We are grateful to the referees for their careful readings and questions, which led to improved clarity of our paper.}

\bigbreak
\section{Notation and Automorphisms} \label {sec:notation}

\RL{We consider only planar embeddings of graphs that are \defn{proper},
which means that every bounded set in the plane intersects only finitely
many vertices and edges.}

We will often work with left-infinite strings of binary digits. We will perform base-2 addition and subtraction with these strings, in which case the last (rightmost) digit is considered to be the units digit, and the place value of each other digit is twice as much as the digit to its right. For example, $$(\dots0011) + 1 = (\dots0100).$$ We also have $$(\dots1111) + 1 = (\dots0000).$$ We will use $\conc$ to denote appending digits to a left-infinite string, for instance, $$(\dots00) \conc 10 = (\dots0010).$$
Thus, we identify left-infinite strings of bits with the dyadic integers, $\Z_2$.
We may also write $a = (a_k)_{k \le 0} = \sum_{k \le 0} a_k 2^{-k}$.
Thus, for example, $a \conc 0 = 2a$.

The symbol $0^k$ will indicate $k$-fold repetition, so $0^3$ is the string $000$.

We are interested in random walks on the following state space. %

\begin{Definition}\label{def:gamma}
Form a disconnected graph $\states$ on the uncountable vertex set $\Z \times \Z_2$ by adding
an edge between the vertices $(m,b)$ and $(n,c)$ if either
\begin{itemize}
\item $m = n$ and $c = b \pm 1$ (such edges are called \defn{horizontal}), or
\item $n = m+1$ and $c = b \conc 0$ or $m = n+1$ and $b = c \conc 0$ (such edges are called \defn{vertical}). 
\end{itemize}
The connected component of the vertex $(0, a)$, with root $(0, a)$, is denoted\/ $\Gamma_{a}$; we endow $\Gamma_a$ with a planar embedding under which the sequence of edges $\bigl((0, a), (0, a-1)\bigr)$, $\bigl((0, a), (1, a \conc 0)\bigr)$, $\bigl((0, a), (0, a+1)\bigr)$ has a positive orientation. While the set $\Z \times \Z_2$ is uncountable, the connected component $\Gamma_a$ is countable, because each vertex has finite degree.
Because every graph $\Gamma_a$ is $3$-connected, every such planar embedding of\/ $\Gamma_a$ is unique up to orientation-preserving homeomorphisms of the plane; see \cite{Imrich} or, for a simpler proof in our context, \cite[Lemma 8.42 and Corollary 8.44]{LyonsPeres}.
The collection of such rooted embedded graphs is denoted\/ $\rstates$; they are in natural bijection with the dyadic integers, $\Z_2$.

We say that the \defn{depth} \RL{or \defn{level}} of the vertex $(m,b)$ is the integer $m$. 
\end{Definition}

\begin{figure}[ht]
\begin{center}
\includegraphics[width=.7\textwidth]{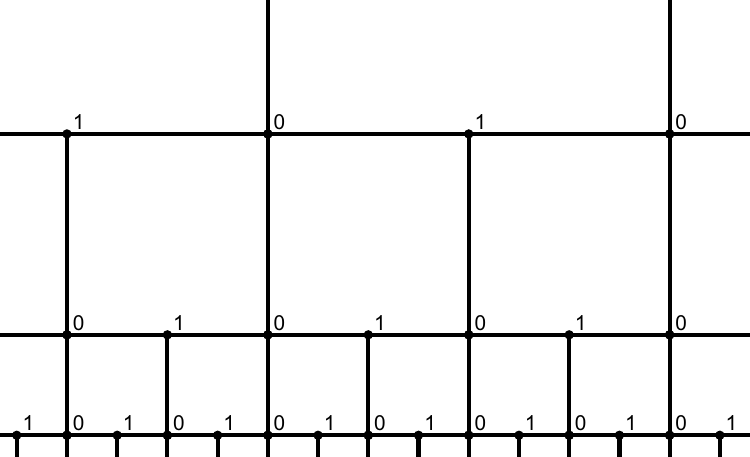}
\end{center}
\caption{The local structure of the graph $\states$. Vertices $\bigl(m, (b_k)_{k \le 0}\bigr)$ are labelled by the last bit $b_0$ of the corresponding left-infinite string, with the entire string $b$ recoverable as in Proposition \ref{prop:readsequence}. Depth $m$ is not notated, but each horizontal slice has depth one greater than the slice above it.} 
\label{fig:gammazoom}
\end{figure}

\begin{figure}[ht]
\begin{center}
\includegraphics[width=.6\textwidth]{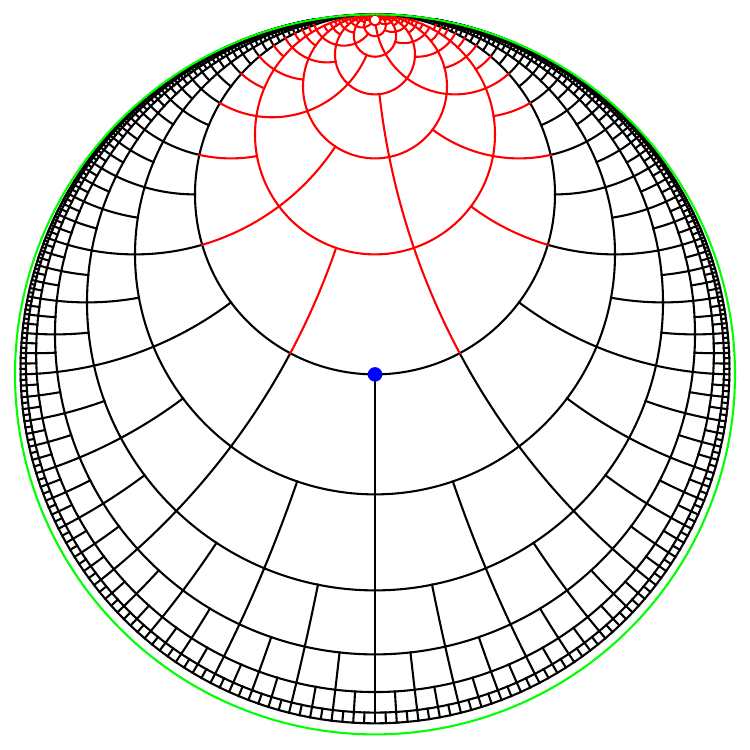}
\end{center}
\caption{A graph $\Gamma_{a}$ in the Poincar\'e disc. The root vertex is coloured blue, and the portion at negative depth is red. The black portion is $\Gamma_+$. In this instance, $a = \dots 110111$. These digits determine the structure of the graph above the root vertex, as described in Proposition \ref{prop:readsequence}.}
\label{fig:gammadisc}
\end{figure}

Figure \ref{fig:gammazoom} shows the local structure of the graph $\states$. %
The global structure of the graphs $\Gamma_{a}$ may be better seen in the Poincar\'e disc, rather than the upper halfplane, as in Figure \ref{fig:gammadisc}. 

For every $a \in \Z_2$, there is a rooted isomorphism from $\Gamma_a$ to $\Gamma_{-a}$ given by $(m, b) \mapsto (m, -b)$. Here, negation is done in the additive group of dyadic integers, $\Z_2$. This isomorphism clearly reverses the orientation specified by Definition \ref{def:gamma}.

\procl r.connected
Two elements $\bigl(m, (b_k)_{k \le 0}\bigr)$ and $\bigl(n, (c_j)_{j \le 0}\bigr)$ of\/ $\Z \times \Z_2$ belong to the same connected component of\/ $\states$ iff there exists $r \le 0$ with $(b_{r-m+\ell})_{\ell \le 0} = (c_{r-n+\ell})_{\ell \le 0}$ or if both $(b_k)$ and $(c_j)$ are eventually 0 or 1, i.e., represent ordinary integers.
\endprocl

A simple random walk $(X_n)_{n \ge 0}$ on $\states$ is defined by simple random walk on the component of the starting vertex $X_0$, that is, each vertex $X_{n+1}$ is a uniformly random neighbour of the preceding vertex, $X_n$. This induces a walk on $\Z_2$ by projecting $(m, b) \mapsto b$; we regard this as a walk on rooted graphs, where we move the root from its present position to one of its neighbours, chosen uniformly at random.
As we will see, there is an orientation-preserving rooted isomorphism between rooted graphs $\Gamma_a$ and $\Gamma_b$ iff $a =b$.

\begin{Definition}\label{def:gammawalk}
The simple random walk on the space $\rstates$ of Definition \ref{def:gamma} moves from the state $\Gamma_{a}$ to the state $\Gamma_{b}$, where $b$ is obtained by choosing uniformly at random from either the following three or four options, depending on whether the last is valid.
\begin{itemize}
\item Adding $1$ to $a$ (referred to as \defn{moving right}).
\item Subtracting $1$ from $a$ (referred to as \defn{moving left}).
\item Appending $0$ to $a$ (referred to as \defn{moving down}).
\item Removing a terminal $0$ from $a$ (only possible if $a_0$ = 0) (referred to as \defn{moving up}).
\end{itemize}
\end{Definition}

The behaviour of this random walk will be analysed by considering the following graph.

\begin{Definition}\label{def:gamma0}
Let\/ $\w$ be the graph whose vertex set consists of all finite strings of binary digits, including the empty string $\e$, with edges between vertices $a$ and $b$ if
\begin{itemize}
\item $b = a \pm 1$ (possible only if $a \ne \e$), or
\item $b = a \conc 0$ or $a = b \conc 0$. 
\end{itemize}
Here, addition of $\pm 1$ to a string of length $n$ is reduced modulo $2^n$, so that for instance $$(111) + 1 = (000).$$
The \defn{depth} \RL{or \defn{level}} of a vertex is the length of its string. 
\end{Definition}

We will mostly be interested in what happens as the depth increases, because the depth of the simple random walk on $\w$ has positive drift (\rref p.driftdown/). Again, $3$-connectivity implies that there are at most two planar embeddings of $\w$ up to planar homeomorphisms, but, in fact, there is only one, because reversing the orientation leads to \RL{an equivalent embedding}: compare Proposition \ref{prop:reflectionauto}.

In comparing the definitions of $\Gamma_{a}$ and $\w$, observe that the depth of a vertex in $\w$ is the length of the corresponding string, which in the finite setting can be read from the string and need not be given as an additional parameter as in the definition of $\Gamma_{a}$.

The relation between the graphs $\Gamma_{a}$ and $\w$ is that a portion of $\Gamma_{a}$ may be `wrapped' to produce $\w$, as follows.

\begin{Remark}\label{rem:wrap}
The portion of\/ $\Gamma_{a}$ induced by all vertices of nonnegative depth is denoted $\Gamma_+$; it does not depend on the choice of $a$ up to (orientation-preserving) rooted isomorphism. If we identify all pairs of vertices $\bigl(m,(b_k)_{k \le 0}\bigr)$ and $\bigl(m,(c_k)_{k \le 0}\bigr)$ of\/ $\Gamma_+$ such that $m \ge 0$ and $b_k = c_k$ for all $k \in (-m, 0]$, then the resulting graph is isomorphic to $\w$.
\RL{We may also express the relationship between $\Gamma_+$ and\/ $\w$ as follows. For each $n \in \Z$, there is an automorphism of\/ $\Gamma_+$ given by $(m, b) \mapsto (m, b + 2^{-m} n)$, where $(m, b) \in \Z \times \Z_2$. The orbit of each vertex in $\Gamma_+$ corresponds to a single vertex of\/ $\w$. We thus refer to $\w$ as the quotient of\/ $\Gamma_+$ by this action of\/ $\Z$.} 
\end{Remark}

This construction of $\w$ does two things to $\Gamma_{a}$---it removes the portion of the graph with negative depth, and then wraps the remaining graph around a cylinder, so that moving sufficiently far to the right or left results in returning to where one started.

\begin{figure}[ht]
\begin{center}
\includegraphics[width=.6\textwidth]{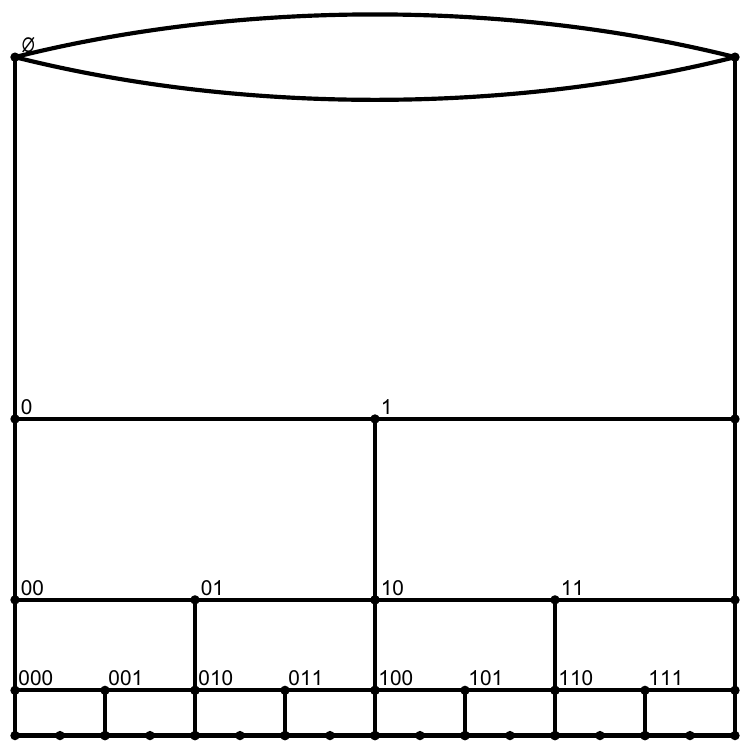}
\end{center}
\caption{The first few levels of the graph $\w$. The left and right boundaries are identified with one another.} 
\label{fig:gamma0}
\end{figure}

A portion of $\w$ is shown in Figure \ref{fig:gamma0} and also in \rref f.wrapped-circ/. Some key details are the following.
\begin{itemize}
\item Each vertex has either degree $3$ or degree $4$. (For the vertex $\e$, we consider each loop as contributing only $1$ to the degree.)
\item From any vertex, it is possible to move to the left or right. Moving in either direction alternates between vertices of degree $3$ and vertices of degree $4$, except from $\e$.
\item From any vertex, it is possible to move down, to a vertex of degree $4$.
\item From a vertex of degree $4$ only, it is possible to move up. This may result in a vertex of degree either $3$ or $4$. 
\end{itemize}

\begin{figure}[ht]
\begin{center}
\includegraphics[width=.6\textwidth]{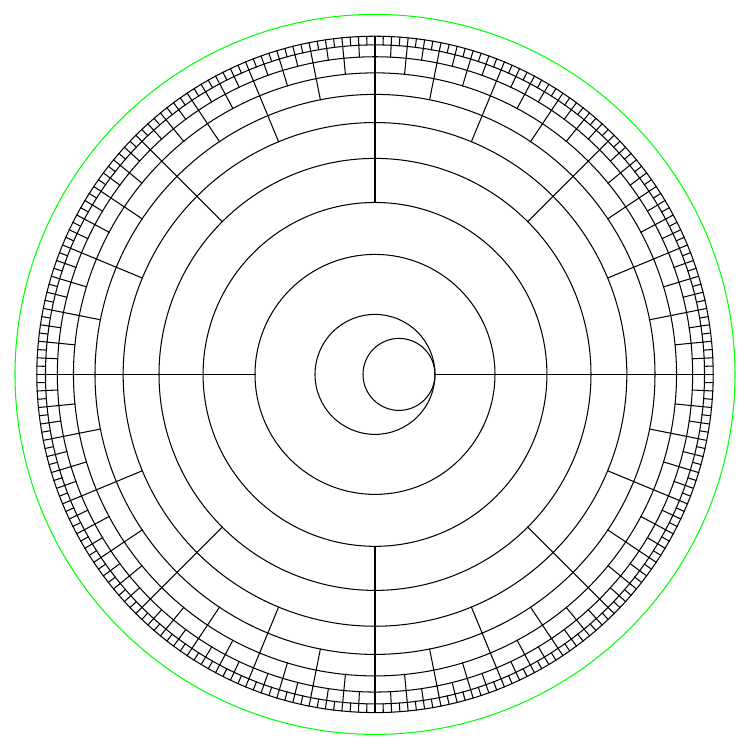}
\end{center}
\caption{The first few levels of the graph $\w$ in a circular drawing in the Poincar\'e disc.}
\label{f.wrapped-circ}
\end{figure}

In the notation of Definition \ref{def:gamma0}, vertices of degree $4$ are those whose strings end with a zero. Moving to the right or left corresponds to adding or subtracting $1$ from the string, moving downward corresponds to appending a zero, and moving upwards to removing a zero from the end of the string. This last operation is only possible when the string ends with a zero, which is why these strings correspond to vertices of degree $4$. (From the empty string, $\e$, the only allowable moves are to stay  at $\e$---in either of two ways---or to append a $0$.)

While the graphs $\Gamma_{a}$ and $\w$ are very structured, they usually do not have many automorphisms.
Our experience is that this comes as a surprise to all who initially hear of it; probably this is because people are used to thinking about $\Gamma_+$, which has automorphism group $\Z \times (\Z/2\Z)$.

\begin{Proposition}\label{prop:reflectionauto}
The graph $\w$ has only one nontrivial automorphism. 
\end{Proposition}
\begin{proof} %
The graph $\w$ has an automorphism $\phi$ that takes a binary string $w$ of length $n$ to the binary string $2^n - w$ of equal length. This may be seen as a reflection in Figure \ref{fig:gamma0}, where it fixes each string of the form $000\dots0$ and $100\dots0$, and interchanges the two horizontal intervals between these fixed points. Equivalently, moves to the left and moves to the right are swapped. 
This is also a reflection in the real axis in \rref f.wrapped-circ/.

This notion of interchanging the notions of right and left is evident in the action of $\phi$ at any fixed depth. The depth $n$ portion of $\w$ may be seen as the Cayley graph of the additive group of integers modulo $2^n$ with generators $\pm1$. There is a unique nontrivial automorphism $\phi_n$ of this graph that preserves the identity, defined by writing an arbitrary element $w$ as a sum of generators and reversing the sign of each of these generators, producing $\phi_n(w) = -w$. This operation of reversing the sign exchanges the roles of the generators $+1$ and $-1$, or equivalently of right and left. 

Any automorphism of $\w$ must fix $\e$, which is the sole vertex with a double loop, must fix $0$, which is the only vertex connected to $\e$, and must fix $1$, which is the only vertex connected to $0$ by two edges. Finally, it must fix the only other neighbours of $0$ and $1$, which are $00$ and $10$. There are only two vertices connected to both $00$ and $10$, namely, $01$ and $11$. Therefore any automorphism of $\w$ either fixes these two vertices or interchanges them.

To see that $\w$ has no nontrivial automorphisms other than $\phi$, it suffices to show that an automorphism of $\w$ that fixes the vertices $01$ and $10$ must fix the entire graph, as automorphisms that swap $01$ and $11$ may be composed with $\phi$. 

If an automorphism of $\w$ fixes each $k$-bit string for some $k \geq 2$, then it fixes each $(k+1)$-bit string that ends with a zero, and then also fixes each $(k+1)$-bit string that ends with a one, because each string ending with a one is adjacent to different pairs of strings ending with zeros, as long as $k \geq 2$. This completes the proof.
\end{proof}

\begin{figure}[ht]
\begin{center}
\includegraphics[width=.6\textwidth]{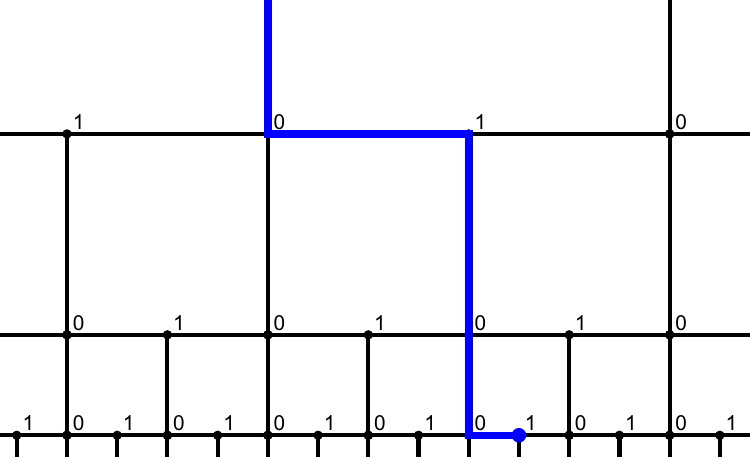}
\end{center}
\caption{A portion of a graph $\Gamma_{a}$ with the path, marked in blue, traced out by the procedure described in Proposition \ref{prop:readsequence}. In this instance, $a = \dots 101$. This string of digits is comprised of the digits labelling only the first vertex at each level along the indicated path.}
\label{fig:stationary}
\end{figure}

Most graphs $\Gamma_{a}$ do not have nontrivial automorphisms, even considering automorphisms of graphs rather than rooted graphs, which need not fix the root vertex.
We will show that edges of $\Gamma_a$ can be classified as either vertical or horizontal from the isomorphism class of $\Gamma_a$, as well as whether traversing a vertical edge moves up or down; whether traversing a horizontal edge moves right or left can then be determined from the planar orientation of the embedding. This will result in the following:

\begin{Proposition}\label{prop:readsequence}
The sequence $(a_i)$ may be read from the oriented graph $\Gamma_a$ via the following procedure. %

\begin{enumerate}
\item Start at the root, and initialise a counter $i$ to be $0$.
\item\label{step:start} If the present vertex has degree $4$, set $a_i$ to be $0$ and step upwards. Otherwise set $a_i$ to be $1$ and step left and then up.
\item\label{step:end} Decrease $i$ by $1$ and repeat steps \ref{step:start} and \ref{step:end} indefinitely. 
\end{enumerate}
The path followed by this procedure in an example is shown in Figure \ref{fig:stationary}.
\end{Proposition}

\begin{Lemma}\label{lem:readhorv}
In any graph $\Gamma_{a}$, the classification of edges as horizontal or vertical may be read from the graph structure, in other words, is preserved by automorphisms.
\end{Lemma}
\begin{proof}
Any horizontal edge is between a vertex of degree $3$ and a vertex of degree $4$, so an edge with two degree-$4$ endpoints must be vertical. 

Consider an edge $e$ between a vertex $x$ of degree $3$ and a vertex $y$ of degree $4$. If $y$ is adjacent to two vertices of degree $4$, then $e$ is a horizontal edge. Otherwise, let $z$ be the unique vertex of degree $4$ adjacent to $y$. If \RL{a} shortest path between $x$ and $z$ that does not go through $y$ has length $3$, then $e$ is horizontal. Otherwise \RL{such a} path has length $7$ and $e$ is vertical. Examples of \RL{these paths are shown in red and green, respectively, in \rref f.horiz-vert/.}
\end{proof}

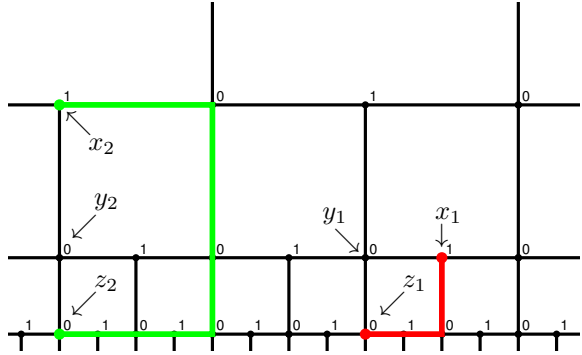
\begin{figure}[ht]
\begin{center}
\begin{overpic}[width=.6\textwidth,tics=10]{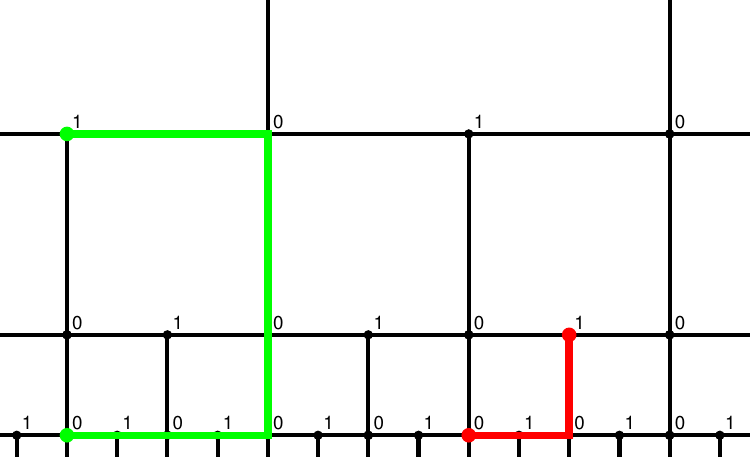}
   \put(74.7,23){$x_1$}
   \put(74.7,19){$\downarrow$}
   \put(55,23){$y_1$}
   \put(57,18){$\searrow$}
   \put(69,11){$z_1$}
   \put(65,7){$\swarrow$}
   \put(14,35){$x_2$}
   \put(10,38.5){$\nwarrow$}
   \put(15,25.5){$y_2$}
   \put(11,20.5){$\swarrow$}
   \put(15,11){$z_2$}
   \put(11,7){$\swarrow$}
\end{overpic}
\end{center}
\caption{A portion of a graph $\Gamma_{a}$ illustrating the proof of Lemma~\ref{lem:readhorv}. \RL{The shortest path from $x_1$ to $z_1$ not passing via $y_1$ has length $3$ and is drawn in red; a shortest path from $x_2$ to $z_2$ not passing via $y_2$ has length $7$ and is drawn in green.}}
\label{f.horiz-vert}
\end{figure}

\begin{Lemma}\label{lem:readdown}
In any graph $\Gamma_{a}$, the notions of `up' and `down' may be read from the graph structure.
\end{Lemma}
\begin{proof}
Let $e$ be a vertical edge between two vertices $x$ and $y$. There are exactly two paths between these two vertices that take one horizontal step, then one vertical step, then two horizontal steps. The initial vertex of these paths is above the final vertex.
\end{proof}

Lemmas \ref{lem:readhorv} and \ref{lem:readdown} imply Proposition \ref{prop:readsequence}, whence the string $a$ is determined from the isomorphism class of the (rooted, oriented) graph $\Gamma_a$. More precisely, there is an orientation-preserving rooted isomorphism between $\Gamma_{a}$ and $\Gamma_{b}$ if and only if $a = b$, as claimed earlier.

\begin{Proposition} \label{prop:evn-per}
A graph $\Gamma_{a}$ has nontrivial automorphisms if and only if the sequence $a$ is eventually periodic.
\RL{In particular, if $a$ is not eventually periodic, then the map $(m, b) \mapsto b$ from vertices of\/ $\Gamma_a$ to $\Z_2$ is injective.}
\end{Proposition}
\begin{proof}
If there is an automorphism that takes $(m, b)$ to $(n, c)$ with $m \ne n$ or $b \ne c$, then by Proposition \ref{prop:readsequence}, either $b = c$ or $b = -c$ (as dyadic integers). Since $(m, b)$ and $(n, c)$ belong to the same connected component of $\states$, it follows from \rref r.connected/ that $b$ is eventually periodic, whence so is $a$.

Conversely, if $a$ is eventually periodic, then $\Gamma_a$ has a vertex $(m, b)$ where $b$ is periodic. Thus, without loss of generality, assume that $a$ itself is periodic. Choose $m < 0$ so that $(a_{m+k})_{k \le 0} = a$. Then there is an orientation-preserving automorphism that takes $(0, a)$ to $(m, a)$, namely, $(n, b) \mapsto (n+m, b)$. 
\Qed

Note that any automorphism of $\Gamma_{a}$ acts on the depths of vertices by addition of a constant, so automorphisms that do not fix the depth of every vertex must have infinite order. Note also that $\Gamma_{a}$ only has a nontrivial depth-preserving automorphism when $\Gamma_{a}$ contains vertices labelled by strings of only zeros (the vertices on the axis of reflection), in which case there are also automorphisms that do not preserve depth. Therefore, the automorphism group of $\Gamma_{a}$ is only ever trivial or isomorphic to $\Z$ or $\Z \times (\Z / 2\Z)$, with the last case arising only when $\Gamma_{a}$ contains vertices labelled by strings of only zeros. 
In no case is $\Gamma_a$ even close to vertex-transitive.

\begin{Remark}\label{rem:bs}
All our graphs $\Gamma_c$ ($c \in \Z_2$) are subgraphs of the standard Cayley graph of the
Baumslag--Solitar group $\BS(1,2) = \pres{a,b}{bab^{-1}=a^2}$. A portion of this Cayley graph is shown in Figure \ref{fig:bscayley}. 

For comparison with our graphs, this Cayley graph is drawn with edges corresponding to the generator $a$ horizontal and edges corresponding to the generator $b$ vertical. The faces in our graphs have degree $5$ and correspond to the relation $bab^{-1}a^{-2} = 1$. This may be understood as saying that from any point, stepping up, right, down, left, and left again results in a cycle, returning to the starting vertex.  

The primary difference is that in $\Gamma_c$, one is permitted to step upwards only from every second
vertex \RL{in a given level}, whereas in the graph of\/ $\BS(1,2)$, it is possible to step upwards from any vertex, but doing so from even or odd vertices result in different branches of the graph. 

A choice of rooted graph $\Gamma_{c}$, identifying the root of\/ $\Gamma_c$ with the identity element
of\/ $\BS(1, 2)$, amounts to choosing a `sheet' in the Cayley graph of the
group $\BS(1,2)$, where it is always possibly to move downward or sideways, but at each level only one of the two upward branches (even or odd) is chosen.
\end{Remark}

\begin{figure}[ht]
\begin{center}
\includegraphics[width=.4\textwidth]{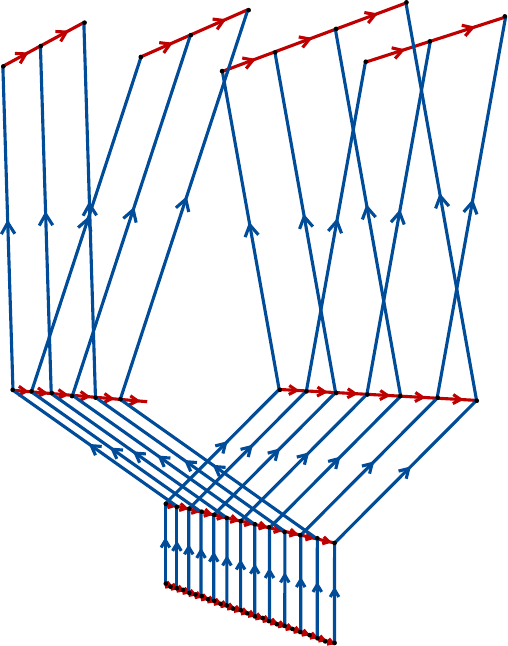}
\end{center}
\caption{A portion of the Baumslag--Solitar group $\BS(1, 2)$ with generator-$a$ edges coloured red
and generator-$b$ edges coloured blue.}
\label{fig:bscayley}
\end{figure}

\bsection{Stationary and Harmonic Measures: Basic Properties}{s.basic}

\RL{We will consider simple random walks on $\states$ and $\w$. Both of these random walks tend downward. The depth of the walk on $\states$ tends to $+\infty$; consequently, on $\w$, the strings that represent the location of the walk have length tending to $\infty$.}

\begin{Definition}\label{def:leaving}
Consider a random walk on either $\states$ or $\w$. For each
integer $n$, let the random variable $T_n$ be the last time for which the
depth of the walk is at most $n$. That is, at time $T_n$, the depth is $n$
and at each time after $T_n$, the depth is greater than $n$. If there are arbitrarily large times for which the depth is at most $n$, then $T_n := \infty$.
We call $T_n$ the \defn{leaving time} for depth $n$.
\end{Definition}

\begin{Proposition}\label{prop:drift}
The simple random walks on $\states$ or $\w$ starting at depth $D_0$ have the property that $(T_{n+1}
- T_n)_{n \ge D_0}$ are IID with distribution not depending on $D_0$, and are independent of $T_{D_0}$. \vadjust{\kern2pt}%
Also, there exists $\lambda > 0$ such that for every $n \ge D_0$, we have $\Ebig{\exp(\lambda (T_{n+1} - T_{n}))} + \Ebig{\exp(\lambda T_{D_0})} < \infty$. 
\vadjust{\kern2pt}%
The distribution of $(T_{n+1}
- T_n)_{n \ge D_0}$ is the same for simple random walks on $\states$ as on $\w$.
\end{Proposition}

To prove this, we first obtain bounds on the drift. From a vertex of degree $4$, it is equally likely that the walk moves up or down, while from a vertex of degree $3$, only down is an option. Thus, the drift comes from moves downward from vertices of degree $3$. 
Write $\depth(x)$ for the depth of a vertex $x$ of \RL{$\states$ or} $\w$.

\procl p.driftdown
Let $(X_n)_{n \ge 0}$ be simple random walk on $\states$ or $\w$, with any choice of $X_0$.
Let $D_t := \depth(X_{2t})$ be the depth of $X_{2t}$. Then
\[
\liminf_{t \to\infty} \frac{D_t}{t} \ge \frac{1}6
\]
and for all $t \ge 12$,
\[
\Pbig{D_t \le D_0 + 1}
\le
e^{-t/1152}.
\]
\endprocl

\rproof
By considering all possible sequences of two steps from the vertex $X_{2t}$, we calculate that
\[
\E[D_{t+1} - D_t \mid X_{2t}] \ge 
\rcases
{1/6 &\text{if $X_{2t}$ has degree 4,}\\
1/3 &\text{if $X_{2t}$ has degree 3.}
}
\]
Thus, $\E[D_{t+1} - D_t \mid X_{2t}] \ge 1/6$ in all cases.
Since the random variables $Z_t := D_{t+1} - D_{t} - \E[D_{t+1} - D_{t} \mid X_{2t}]$ are uncorrelated and take values in $[-2, 2]$ for $t \ge 0$, the strong law of large numbers yields that their average tends to 0 a.s.\RL{; see, e.g., \cite[Theorem 13.1]{LyonsPeres}}. This gives the first result. Moreover, $Z_t$ are martingale differences, whence the Hoeffding--Azuma inequality yields
\[
\P[D_t \le D_0 + 1]
\;\le\;
\PBig{\sum_{m=0}^{t-1} Z_m \le 1-t/6}
\;\le\;
\PBig{\sum_{m=0}^{t-1} Z_m \le -t/12}
\;\le\;
e^{-t/1152}.
\qedhere
\]
\Qed

\begin{proof}[Proof of Proposition \ref{prop:drift}]
Note that each level is reached for the first time at a vertex of degree 4, except possibly for the starting level.
The portion $\Gamma_+$ of the graph at the levels equal to or greater than that of a given
vertex and rooted at that vertex does not depend on that vertex or on which graph in $\states$ the walk takes place, in the sense
that it is the same up to rooted isomorphism, so $T_{k+1} - T_k$ are IID in $k \ge D_0$ and have distribution independent of $D_0$.
\RL{Likewise, for the random walk on $\w$, the random variables $T_{k+1} - T_k$ are IID in $k \ge D_0$;}
in addition, a walk on $\Gamma_+$ projects to a walk on $\w$ in a way that preserves changes in levels, whence the distribution of $T_{k+1} - T_k$ is the same on $\states$ as on $\w$.
Finally, let $F_n$ be the first time the random walk is at level $n$. For $u \ge 1$,
\[
\P[T_{k+1} - T_k \ge u]
\;\le\;
\P[T_{k+1} - F_{k+1} \ge u - 1]
\;\le\;
\sum_{t \ge u - 1} \P[\depth(X_{t + F_{k+1}}) = \RL{D(X_{F_{k+1}})}]
\]
and
\[
\P[T_{D_0} \ge u]
\;\le\;
\sum_{t \ge u} \P[\depth(X_t) = \RL{D_0}],
\]
whence exponential tail bounds are provided by \rref p.driftdown/.
\end{proof}

We are interested in the following features of the limiting behaviour of the simple random walks on $\rstates$ and $\w$.
\RL{As mentioned earlier,} we sometimes identify $\rstates$ with the set of dyadic integers.

\begin{Definition}
A \defn{stationary (probability) measure} for simple random walk on $\rstates$ is a probability measure on the dyadic integers that is stationary for the induced random walk. We will show that there is a unique such probability measure and denote it by $\sm$.
\end{Definition}

We will show in Proposition \ref{p.hmexists} that simple random walk on $\w$, considered as a sequence of finite binary strings, converges coordinatewise a.s.\ to a right-infinite string. Any given infinite string will be such a limit with probability 0 \RL{by Lemma~\ref{lem:nopointmasses}}, so there is no loss of information if we identify a right-infinite string with the number in $[0, 1)$ of which it is a binary representation.
\RL{With this identification,}
the convergence of a path in $\w$ to an element of $[0,1)$ may be
seen as convergence of the horizontal \RL{position} in Figure
\ref{fig:gamma0}. It will be convenient to use this \RL{position} to describe
certain vertices---for example, the vertices $01 \conc 0^k$ are at
$\sfrac{1}{4}$ for all $k \ge 0$. %
We may also identify $[0, 1)$ with $\R/\Z$.

\begin{Definition}
The \defn{harmonic measure} for simple random walk on $\w$ starting at the vertex $\e$ is the probability measure on the real interval $[0,1)$ that is the law of the limit of the horizontal \RL{positions} of the walk. We use $\hm$ for this measure.
\end{Definition}

To interpret the horizontal \RL{position} of the walk as a point in $[0,1)$, refer to Figure \ref{fig:gamma0}. Horizontal steps at depth $n$ have size $2^{-n}$, and the graph is wrapped around a cylinder so that the left and right sides are identified---these are the horizontal positions $0$ and $1$.

Because $\w$ is the quotient of $\Gamma_+$ by $\Z$, we may denote the vertices in $\Gamma_+$ by $(n, b)$, where $n \in \Z$ and $b$ is a vertex of $\w$. Here, the inverse image of $b$ under the quotient map is $\Z \times \{b\}$, which are assumed ordered from left to right. If $b$ has horizontal \RL{position} $x \in [0, 1)$, then we identify $(n, b)$ with $n + x \in \R$.
\RL{We will also refer to $n + x$ as the horizontal position of $(n, b)$.}

\begin{Definition}
The \defn{harmonic measure} for simple random walk on $\Gamma_+$ starting at the vertex $(0, \e)$ is the probability measure $\hp$ on $\R$ that is the law of the limit of the horizontal
\RL{positions} of the walk.
\end{Definition}

\RL{The quotient map $(n, b) \mapsto b$ from $\Gamma_+$ to $\w$ pushes forward the random walk on $\Gamma_+$ to the random walk on $\w$, whence it also pushes forward $\hp$ to $\hm$.}
We sometimes regard $\hm$ as a measure on $\{0, 1\}^{\N}$ and correspondingly $\hp$ as a measure
on $\Z \times \{0, 1\}^{\N}$.

\begin{Proposition} \label{p.hmexists}
The harmonic measures $\hm$ and $\hp$ are well defined, in the sense that random walks on $\w$ and on $\Gamma_+$ can be considered to converge a.s.\ to elements of\/ $[0,1)$ or $\R$, respectively.
Furthermore, the support of $\hm$ is $\R/\Z$ and for every path $(\e, x_1, x_2, \ldots, x_k)$ in $\w$ and every $\epsilon > 0$, the probability is positive that the random walk on $\w$ starting with $X_0 = \e$ has $X_1 = x_1, \ldots, X_k = x_k$ and has the property that for all $n \ge k$, the horizontal position of $X_n$ differs from the horizontal position of $X_k$ by less than $\epsilon$.
\end{Proposition}
\begin{proof}
Because $\w$ is a quotient of $\Gamma_+$, it suffices to prove that $\hp$ is well defined in order to prove that $\hm$ is well defined.
For $n \ge 0$, let $Z_n$ be the number of times that the random walk is at level $n$ and $H_n$ be the total (signed) change in horizontal \RL{position} while at level $n$. At most $Z_n$ sideways steps are taken on level $n$, and each of these \RL{changes the horizontal position by} $\sfrac{1}{2^n}$, so $|H_n| \le Z_n/2^{n}$. By Proposition \ref{prop:drift}, $s := \E[T_{n+1} - T_n] \vee \E[T_0] < \infty$, so $\E[T_n] \leq s(n+1)$. \RL{Because $Z_n \le T_n$,} it follows that $\Ebig{\sum_n Z_n/2^{n}} \le \sum_n s (n+1)/2^{n} < \infty$. Therefore $\sum_n H_n$ converges a.s., which is the result.
Because the tails of $\sum_n H_n$ are arbitrarily small, the rest of the proposition follows.
\end{proof}

\RL{An almost identical argument shows the following:
\procl l.sidefuzz
For every $\epsilon > 0$, there is an $n \ge 0$ such that for every $a$, the simple random walk $(X_t)_{t \ge 0}$ on $\Gamma_a$ starting at $X_0 = (0, a)$ has the property that with probability at least $1 - \epsilon$, we have for all $t \ge n$ that $X_t \in \Gamma_+$ and the horizontal position of $X_t$ differs from the horizontal position of $X_n$ by less than $\epsilon$. \qed
\endprocl}

\begin{Lemma}\label{lem:nopointmasses}
For every $x \in [0,1)$, $\>\hm(x) = 0$ and
for every $x \in \R$, $\>\hp(x) = 0$.
\end{Lemma}
\begin{proof}
It suffices to prove the second statement.
Fix $x \in \R$.
Define $S_n$ to be the first time after $T_n$ that the walk makes a horizontal step $s_n = \pm 1$, and let $L_n$ be the level at which the walk makes that step.
Because each $T_n$ is a.s.\ finite, we may choose an increasing sequence $(n_k)_{k \ge 1}$ so that for each $k$, we have $\P[S_{n_k} < T_{n_{k+1}}] \ge \hp(x)/2$.
Let $A$ be the set of paths that tend to $x$, so that $\P(A) = \hp(x)$.
Let $A_k$ be the set of paths in $A$ for which $S_{n_k} < T_{n_{k+1}}$, so that $\P(A_k) \ge \P(A)/2$.
Write $A'_k$ for the sequences obtained from $A_k$ by changing the step $s_{n_k}$ to $-s_{n_k}$; all sequences in $A'_k$ still correspond to paths in $\Gamma_+$, and $\P(A'_k) = \P(A_k)$. 
However, paths in $A'_k$ tend to $x - s_{n_k}/2^{L_{n_k}-1}$. Note that this limit depends on $L_{n_k}$, which is not the same for all paths in $A'_k$. However, by the definition of the sequence $(n_k)$, paths in $A'_j$ and $A'_l$ cannot have the same limit for $j \ne l$. Therefore the sets $A'_k$ ($k \ge 1$) have disjoint sets of limits in $\R$, yet each of these sets has probability at least $\hp(x)/2$. Therefore, $\hp(x) = 0$, as desired.
\end{proof}

\procl t.stationary
There is a unique stationary measure $\sm$ on $\rstates$ for simple random walk. \RL{Furthermore, $\sm$ is continuous, that is, every singleton has measure $0$.}
\endprocl

\rproof
The key idea is to consider a random walk on $\Gamma_a$ as developing a longer and longer past history, which converges to the stationary measure. \RL{Every time the random walk leaves a level for the last time, the future moves of the random walk are independent of the graph above that level, and therefore have the same distribution no matter what $a$ is. The existence of such regeneration times shows that there is at most one stationary measure. Furthermore, we can construct a stationary measure from them as follows.} We break up the walk into segments between the last times it leaves successive levels. These segments are IID, and have finite expected length by Proposition \ref{prop:drift}. \RL{The usual length-biasing and uniform choice then yields a stationary process of random walk moves. These random walk moves must still be converted to dyadic integers.} Whether the walk is presently at a vertex of degree $3$ or degree $4$ depends only on the most recent (partial) segment. If it is at a vertex of degree $4$, then the degree of the vertex immediately above depends only on the most recent two segments, and so on. Considering increasing numbers of these segments will allow us to obtain the stationary measure as a limit.   
\RL{Lastly, the IID segments of random walk moves have distribution that is invariant under switching left and right moves; this leads to continuity of $\sm$.}

\RL{We now begin the proof.}
Consider simple random walk $(X_t)_{t \ge 0}$ on $\Gamma_a$ starting at $(0, a)$.
Each transition is a move on $\Gamma_a$ that is either left, right, up, or down from the current position; code this as a symbol $M_t \in \{\mathrm L, \mathrm R, \mathrm U, \mathrm D\}$, where the walk moves from $X_{t-1}$ to $X_t$ via $M_t$.
The sequences $(M_{T_n + t})_{1 \le t \le T_{n+1} - T_n}$ are IID for $n \ge 0$. (Here, $M_{T_n + 1} = \mathrm D$ for all $n$.) Let $\nu$ denote their common law.
Because $\E[T_{n+1} - T_n] < \infty$, we may define $\mu_1$ to be the distribution on $\Z^+$ defined by 
\[
\mu_1(k)
:=
k \P[T_1 - T_0 = k]/\E[T_1 - T_0]
\]
and then $\mu_2$ to be the law of $V$ when $V$ is uniformly chosen from $\{0, 1, 2, \ldots, W-1\}$ and $W$ has law $\mu_1$.
Let $\nu_0$ denote the law of $(M_{T_0 + t})_{1 \le t \le V}$, where $V$ is independent of the walk and has law $\mu_2$.
\RL{Consider now sequences of $\{\mathrm L, \mathrm R, \mathrm U, \mathrm D\}$ indexed by the nonpositive integers.}
Then renewal theory shows that \RL{the finite-dimensional distributions of $(M_{t+s})_{-t \le s \le 0}$ tend as $t \to\infty$ to those of the concatenation of $S_n$ ($n \le 0$)}, where $S_n$ are independent for $n \le 0$ with law $\nu$ when $n < 0$ and with law $\nu_0$ when $n = 0$.
In particular, this is independent of $a$.
\RL{Indeed, the graph $\Gamma_a$ has cycles of length $5$, whence the distribution of $T_1 - T_0$ is nonlattice. Therefore,
the discrete key renewal theorem shows that for the delayed renewal process $(T_n)_{n \ge 0}$ with renewals at $T_n$, the age distribution tends to $\nu_0$. This gives convergence to the equilibrium renewal process; e.g., \cite[Theorem 3.1]{Lindvall}. Usually one looks into the future, but one can just as well look into the past, as we do here. Given this, the fact that the sequences $(M_{T_n + t})_{1 \le t \le T_{n+1} - T_n}$ are IID yields our claim.}

\RL{We next convert these random walk moves to dyadic integers. 
Each} finite sequence $\sigma$ from $\{\mathrm L, \mathrm R, \mathrm U, \mathrm D\}$ corresponds to a function $f_\sigma \colon \Z_2 \to \Z_2$ by
composing the individual-symbol functions $f_{\mathrm L}(b) := b - 1$, $f_{\mathrm R}(b) := b + 1$, $f_{\mathrm U}(b) := b/2$, $f_{\mathrm D}(b) := 2b$,
provided $f_{\mathrm U}$ is applied only to even dyadic integers. \RL{Here, for a sequence $\sigma = (m_1, m_2, \ldots, m_k)$, we compose in the order $f_\sigma := f_{m_k} \circ \cdots \circ f_{m_2} \circ f_{m_1}$.} Say that a finite sequence $\sigma$ is `definable' if \RL{$f_\sigma$} satisfies that restriction
when applied to \RL{every dyadic integer} and is `permissible' if for \RL{every} nonempty initial segment of $\sigma$, the number of symbols $\mathrm D$ is strictly
larger than the number of symbols $\mathrm U$\RL{; in particular, the first symbol of each permissible $\sigma$ is $\mathrm D$}. The sequences $S_n$ above are guaranteed to be both definable and permissible. Furthermore, the value of 
$f_\sigma(b)$ mod 2 does not depend on $b$ for definable and permissible $\sigma$. More generally, $f_{\sigma_1} \circ f_{\sigma_2} \circ \cdots \circ f_{\sigma_k}(b)$ mod
$2^k$ does not depend on $b$ for any definable and permissible sequences $\sigma_1, \ldots, \sigma_k$. Therefore, $f_{S_0} \circ f_{S_{-1}} \circ \cdots \circ f_{S_{n}}(b)$ has a limit a.s.\ as $n \to - \infty$ and does not depend on $b$; its law is $\sm$. 

\RL{Finally, for $n < 0$, the moves in $S_n$ result in a net change to the depth of $1$ and in a net horizontal change compared to the location of $M_{T_n + 1}$. That is, we may write $\bar X_{T_{n+1}} = 2 \bar X_{T_n} + H_n$ for some IID 
\vadjust{\kern2pt}%
$\Z$-valued random variables $H_n$, where $X_t = \bigl(\depth(X_t), \bar X_t\bigr)$ with $\depth(X_t)$ the depth of $X_t$. The law of $\nu$ is invariant under switching $\mathrm L$ and $\mathrm R$, whence $-H_n$ has the same law as $H_n$. The law $\nu^*$ of the limit of $f_{S_{-1}} \circ f_{S_{-2}} \circ \cdots \circ f_{S_{n}}(b)$ is that of $\sum_{n < 0} 2^{1-n} H_n$, a series that converges in the $2$-adic metric by the estimates in the proof of Proposition~\ref{p.hmexists}. Now an argument very similar to that proving Lemma~\ref{lem:nopointmasses} shows that the law of $\nu^*$ is continuous, whence so is $\sm$.}
\Qed

Uniqueness of the stationary measure guarantees that $\sm$ is ergodic, i.e., every measurable set $A
\subseteq \Z_2$ that is closed for the random walk has measure 0 or 1.

\procl c.nonperiodic
The set of eventually periodic $a \in \Z_2$ has $\sm$-measure 0.
\endprocl

\rproof
There are only countably many finite binary strings, \RL{whence there are only countably many eventually periodic strings. Thus, the result is immediate from \rref t.stationary/}.
\Qed

The following result shows that $\hm$ can be seen in the `tail' of $\sm$.
It will be important for our proof of singularity (\rref t.curien/).

\procl p.smhm
Let $\sigma$ be a finite binary string of length $\ell$.
For $\sm$-a.e.\ $a \in \Z_2$, 
\begin{align*}
\lim_{n \to \infty} \frac1n &\bigl|\bigl\{ k \in [0, n-1] : (a_{-k-\ell+1}, a_{-k-\ell+2}, \ldots,
a_{-k-1}, a_{-k}) = \sigma \bigr\}\bigr|
\\ &=
\hm\bigl(\{b \in \{0, 1\}^{\N} : (b_0, b_1, \ldots, b_{\ell-1}) = \sigma\}\bigr)
=: \hm(\sigma).
\end{align*}
\endprocl

\rproof
\RL{\rref c.nonperiodic/ says that almost no $a$ are eventually periodic, so for $\sm$-a.e.\ $a \in \Z_2$, the map $\depth_a(b) := m$ for vertices $(m, b)$ of $\Gamma_a$ is well defined in view of Proposition~\ref{prop:evn-per}.}
Consider the stationary simple random walk $(X_t)_{t \ge 0}$ on $\rstates$.
If $X_0 = a$, then \RL{for each $t \ge 0$, there is some $m \in \Z$ such that $(m, X_t)$ is a vertex of $\Gamma_a$; with the preceding notation, we may write $m = \depth_a(X_t) = \depth_{X_0}(X_t)$. Given $X_0 = a$, we may ($\sm$-a.s.) identify the walk $(X_t)_{t \ge 0}$ with simple random walk on $\Gamma_a$ starting at its root, $(0, a)$.}
Now consider the standard two-sided \RL{stationary} extension $(X_t)_{t \in \Z}$ of the stationary simple random walk on
$\rstates$.
Thus, \RL{$(X_{t + t_0})_{t \ge 0}$} has the same law as $(X_t)_{t \ge 0}$ for every $t_0 \in \Z$. Because the walk \RL{on $\Gamma_{X_0}$} drifts downward, this two-sided extension almost surely \RL{has the property that for every $m \in \Z$, the number of $t \in \Z$ with $\depth_{X_0}(X_t) = m$ is finite.}
Let $(Z_1, Z_2, \ldots, Z_N)$ be the list of times $t \in \Z$ with $\depth_{\RL{X_0}}(X_t) = \depth_{\RL{X_0}}(X_0) = 0$,
listed in increasing order. Thus, $Z_N = T_0$\RL{, with $T_n$ the leaving times of Definition \ref{def:leaving}}.

Write $B_{k}(a) := (a_{-k-\ell+1}, a_{-k-\ell+2}, \ldots, a_{-k-1}, a_{-k})$.
\RL{Lemmas~\ref {Lemma:sidefuzz} and \ref{lem:nopointmasses}, together with stationarity, imply that}
\[
\lim_{k \to\infty} \Pbig{B_{k}(X_{Z_1}) = B_{k}(X_{Z_2}) = \cdots = B_{k}(X_{Z_N})} = 1.
\]
Because one of the times $Z_1, \ldots, Z_N$ is 0, we obtain
\[
\lim_{k \to\infty} \P\bigl([B_{k}(X_0) = \sigma] \xor [B_{k}(X_{Z_N}) = \sigma]\bigr) 
=
0.
\]
Let $A_{k, n}$ be the event that $B_{k}(X_{T_n}) = \sigma$; \RL{using the definition of $A_{k, 0}$ and the fact given earlier that $Z_N = T_0$,} the preceding equation can be
written as
\rlabel e.stopornot
{\lim_{k \to\infty} \P\bigl([B_{k}(X_0) = \sigma] \xor A_{k, 0}\bigr) 
=
0.}
\RL{Let $L$ be the set of random walk trajectories where the level at time $0$ is never visited again. The law of $(X_t)_{t \in \Z}$ is invariant and ergodic under the left shift because $\sm$ is invariant and ergodic, and the leaving times $T_n$ correspond to shifts that bring the trajectory to $L$. Because the return map to $L$ is also measure-preserving and ergodic, it follows that the sequence 
$(\I{A_{k,n}})_{n \in \Z}$ is stationary and ergodic.}
Therefore, for each $k \ge 0$, the ergodic theorem yields
\rlabel e.statCesaro
{\lim_{n \to \infty} \frac1n \bigl|\bigl\{ m \in [0, n-1] :  A_{k, -m} \bigr\}\bigr|
=
\P(A_{k, 0}) \quad \text{a.s.}}
On the other hand, 
Proposition \ref {p.hmexists} shows that
$\I{A_{n-\ell, n}}$ converges a.s.\ as $n \to\infty$ with
$
\lim_{n \to\infty} \P(A_{n-\ell, n}) = \hm(\sigma)
$\RL{; since
by stationarity, $\P(A_{k, n})$ is the same for all $n$, we obtain}
\rlabel e.Ak0sigma
{\lim_{k \to\infty} \P(A_{k, 0}) = \hm(\sigma).}
Since $(\I{A_{n-\ell, n}})_{n \ge 0}$ is a Cauchy sequence a.s., we have
\[
\lim_{k \to\infty} \sup_{m, n \ge k} \P(A_{n-\ell, n} \xor A_{m-\ell, m}) = 0.
\]
By stationarity, $\P(A_{n-\ell, n+r} \xor A_{m-\ell, m+r})$ does not depend on $r$.
Therefore,
\[
\lim_{k \to\infty} \sup_{m, n \ge k} \sup_{r \in \Z} \P(A_{n-\ell, n+r} \xor A_{m-\ell, m+r}) = 0.
\]
Choosing $n = k + \ell$, $\> m=k+\ell+s$, and $r = -k -\ell - s$ yields
\rlabel e.shifted
{\lim_{k \to\infty} \sup_{s \ge 0} \P(A_{k, -s} \xor A_{k + s, 0}) = 0.}
Combining \rref e.statCesaro/, \rref e.Ak0sigma/, and \rref e.shifted/, we obtain
\[
\lim_{k \to\infty}
\lim_{n \to \infty} \frac1n \bigl|\bigl\{ m \in [0, n-1] :  A_{k+m, 0} \bigr\}\bigr|
=
\hm(\sigma) \quad \text{a.s.},
\]
which is the same as
\[
\lim_{n \to \infty} \frac1n \bigl|\bigl\{ m \in [0, n-1] :  A_{m, 0} \bigr\}\bigr|
=
\hm(\sigma) \quad \text{a.s.}
\]
In view of \rref e.stopornot/, we may write this as
\[
\lim_{n \to \infty} \frac1n \bigl|\bigl\{ m \in [0, n-1] :  B_{m}(X_0) = \sigma \bigr\}\bigr|
=
\hm(\sigma) \quad \sm\text{-a.s.}
\]
This is the desired result.
\Qed

Unfortunately, we do not have an explicit description of $\sm$. For instance, one may ask about the proportion of time a simple random walk spends on vertices of degree three.

\begin{Definition}
Given a random walk on $\states$ or $\w$, denote by $p_3$ the limiting fraction of the time spent at vertices of degree $3$.
\end{Definition} 

The limit $p_3$ need not exist for general random walks on these graphs, but it will for the simple random walks we study.

\begin{Proposition}\label{prop:p3exists}
For simple random walks on $\states$ or $\w$, the limit $p_3$
exists a.s., is constant, and is the same for $\states$ as for
$\w$. 
\end{Proposition}
\begin{proof}
Consider a simple random walk on $\states$ or $\w$ that starts at a vertex at
depth $n$. This random walk may be broken up into the segments between the
leaving times $T_k$ and $T_{k+1}$ for $k \geq n$, in addition to the initial
segment up to time $T_n$. By Proposition \ref{prop:drift}, the expected time in
each such interval spent at vertices of degree $3$ is bounded,
and the times spent at vertices of degree $3$ in each of these intervals are IID, except for the initial segment. Therefore the limit $p_3$ exists and is equal to the ratio of the average time spent at vertices at degree $3$ between times $T_k$ and $T_{k+1}$ to the average interval length $T_{k+1} - T_k$.
\end{proof}

We may also write $p_3 = \sm\bigl\{(a_k)_{k \le 0} : a_0 = 1\bigr\}$,
which equals $\P[\deg X_n = 3]$ for every $n \ge 0$ if $(X_n)_{n \ge 0}$ is the random walk on $\Gamma_a$ when $a$ has law $\sm$.
This is a basic quantity in
understanding the behaviour of simple random walk on $\states$.
For example, the speed (drift downwards) of simple random walk on every graph in $\rstates$ or on $\w$
equals $p_3/3$ because the drift at each vertex of degree 3 is 1/3 while the drift at each vertex of
degree 4 is 0.
We can get some simple bounds on $p_3$ as follows.

Each vertex of $\Gamma_{a}$ either has degree $3$ or degree $4$. If it has degree $4$, then the
vertex above it has either degree $3$ or $4$. Let $p_{4,3}$ and $p_{4,4}$ be the proportions of time
spent at vertices of degree $4$ whose upper neighbour has the appropriate degree. These quantities
exist by similar arguments to Proposition \ref{prop:p3exists}. Note that 
\rlabel e.p4s
{p_3 + p_{4,3} + p_{4,4} = 1.}

Consider a vertex chosen according to the stationary measure, and take a single step. The resulting
measure is still the stationary measure, so considering the probability of being at a vertex of
degree $3$ gives us that 
\begin{align} \label{e.p43}
p_3 &= p_3 \cdot 0 + p_{4,3} \cdot \tfrac34 + p_{4,4} \cdot
\tfrac12 \\ & \leq (1-p_3)\cdot \tfrac34. \nonumber
\end{align}
Therefore $p_3 \leq \sfrac{3}{7}$.
It follows that the downward speed is at most 1/7: 

\procl p.speed
The drift downwards of simple random walk on every graph in $\rstates$ or on $\w$ is at most $1/7$.
\qed
\endprocl

\begin{Remark} \label{r.moreonspeed}
The same calculation shows that 
\[
p_3 \geq (1-p_3) \cdot \tfrac12,
\]
so $p_3 \geq \sfrac{1}{3}$ and the speed is at least $1/9$ (the bound of \rref p.driftdown/ is $1/12$).

Using the equations \eqref {e.p4s} and \eqref{e.p43}, we find that
\[
p_{4, 3} = 6 p_3 - 2 \quad \text{ and } \quad p_{4, 4} = 3 - 7 p_3.
\]
Similarly, with $p_{4,4,3}$ being the probability of being at a degree $4$ vertex with a degree $4$ vertex above and a degree $3$ vertex above that, we have 
\[
p_{4,3} = p_{4,4,3}\cdot\tfrac14 + p_3\cdot\tfrac23.
\]
Substituting the previous expression for $p_{4,3}$ gives that $p_{4, 4, 3} = \frac{64}3 p_3 - 8$, which
implies the slightly better bounds $p_3 \in (3/8, 27/67)$ since $0 < p_{4, 4, 3} < 1 - p_3$. 
One might hope to perform increasingly more detailed versions of this calculation, obtaining better and better bounds. This does not seem feasible because when vertices are classified by their neighbourhoods of radius $r$ in this way, the number of different types of vertex grows exponentially in $r$. For instance, our first calculation used the fact that taking a horizontal step from a vertex of degree $3$ results in a vertex of degree $4$, and the next that taking a horizontal step from a vertex of degree $3$ is equally likely to contribute to $p_{4,3}$ or $p_{4,4}$. However, at the next level of detail, we would need to consider two different types of vertices of degree $3$---those where a sideways step contributes to $p_{4,4,3}$, and those where it contributes to $p_{4,4,4}$.

In Section \ref{sec:numerics}, we will give much more precise numerical estimates by other methods.

\end{Remark}

We now compare random walks on $\states$ and $\w$ to random walks
on the dual graphs, where the `dual graph' of $\states$ is taken as the union of the dual graphs of
the connected components of $\states$. The general shape of these dual graphs was shown in Figure
\ref{fig:gammadual-intro}. The appropriate analogue of $\w$ in the dual setting is not the plane dual of $\w$, but rather is obtained from the dual of $\Gamma_a$ in the same way that $\w$ was obtained from $\Gamma_a$. Consider the dual of any graph $\Gamma_a$, and define the graph $\wdd$ to be the portion below any vertex, then `wrapped' mod 1. This is shown in \rref f.gamma0dual/.
To avoid ambiguity, we will sometimes refer to $\states$, $\Gamma_+$, and $\w$ as the \defn{primal} graphs, in contrast with their dual graphs.

\begin{figure}[ht]
\begin{center}
\includegraphics[width=.4\textwidth]{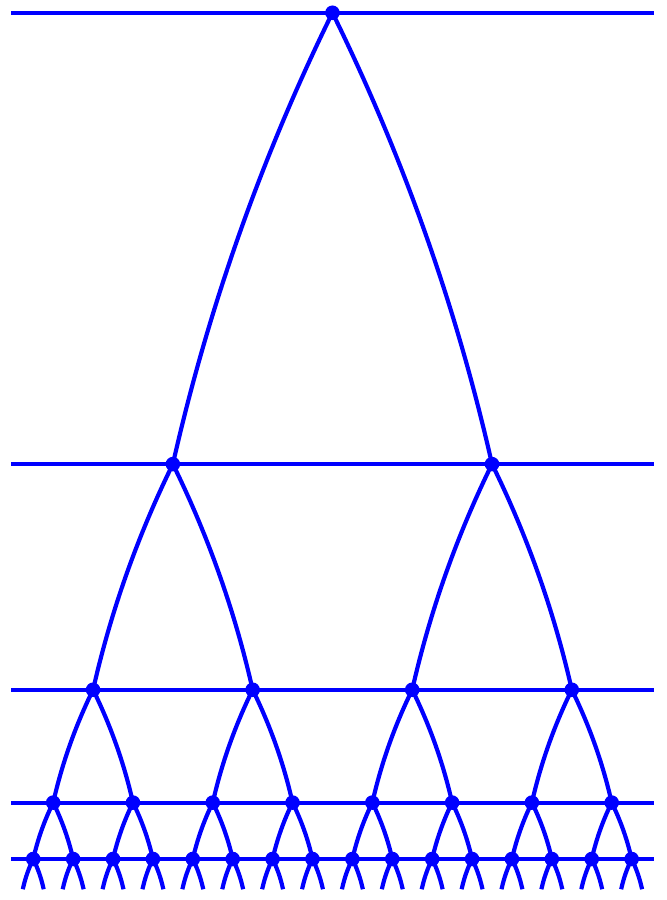}
\end{center}
\caption{A portion of the graph $\wdd$, the dual of any of the graphs $\Gamma_{a}$ below any vertex and wrapped. Here, the edges going off the left side come back in on the right side, except that the top edge is a double edge.}
\label{f.gamma0dual}
\end{figure}

In $\states$, the faces are rectangular, each with five edges---one edge each on the top, left, and right sides, and two edges below. Therefore a step of the dual walk can go up, left, right, down-left or down-right.
Denote by
$\rstatesd$ the collection of oriented, rooted, plane graphs dual to $\Gamma_a \in \rstates$.
It will be convenient to label the faces of $\states$ (and hence the vertices of $\rstatesd$ and of
$\wdd$) with binary strings (dyadic integers) in the same way as the vertices. We will label each face with the string labelling the vertex in its upper left corner, as shown in Figure \ref{fig:gammaduallabels}.

\begin{figure}[ht]
\begin{center}
\includegraphics[width=.6\textwidth]{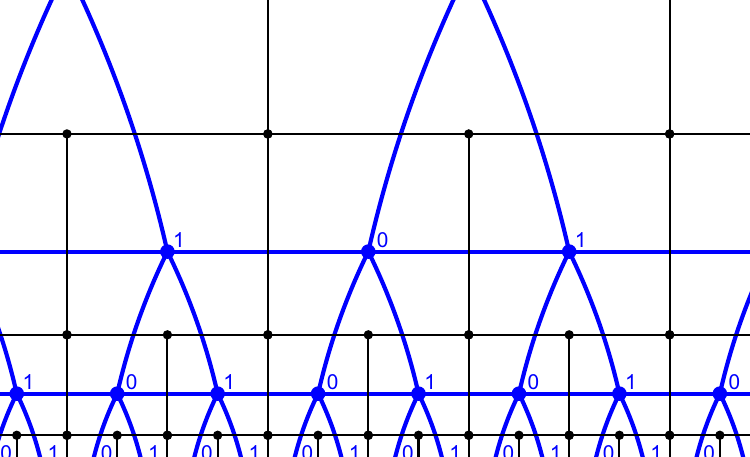}
\end{center}
\caption{The vertex labels of the dual graph. These labels are the same as those of the primal vertex above and to the left of each dual vertex.}
\label{fig:gammaduallabels}
\end{figure}

With these labels, the dual walk on binary strings is defined as follows.

\begin{Definition}
To take a step of the dual walk, perform one of the following five operations, chosen uniformly at random.
\begin{itemize}
\item Add or subtract $1$.
\item Remove the final bit.
\item Append a $0$ or a $1$.
\end{itemize}
\end{Definition}

Notice that the dual walk always allows a step upwards (removing the final bit), and has two different ways to step downward (appending either a $0$ or a $1$)---compare to the primal walk of Definition \ref{def:gammawalk}.
It is clear that the dual walk drifts downwards at speed 1/5.

We will also consider the stationary measure and the harmonic measure for the dual walk. 

\begin{Definition}
A \defn{stationary (probability) measure} for simple random walk on $\rstatesd$ is a probability measure on the dyadic integers that is stationary for the induced random walk. 
Arguments similar to those used to prove \rref t.stationary/ show that there is a unique stationary
measure, $\smd$.
\end{Definition}

In contrast to the primal case, we can easily identify $\smd$.
Namely, it is the symmetric Bernoulli process on $\{0, 1\}^{-\N}$, i.e., bits are independent fair coin flips.

\procl p.sdmuniform
The measure $\smd$ on $\{0, 1\}^{-\N}$ is the product measure $\bigl(\frac{\delta_0 + \delta_1}2\bigr)^{\otimes -\N}$.
\endprocl

\rproof
This product measure is the same as Haar measure on the
compact group, $\Z_2$. Being Haar
measure, it is invariant under both adding 1 and under subtracting 1, which correspond to the walk moving right or left.
This measure is also clearly invariant under deleting the last bit, which happens when
the walk moves up, and under concatenating a random fair bit, which is
the fair mixture of moving down-left or down-right. 
Since this measure is invariant under each of these transformations, it is also invariant under the mixture of them given by a step of simple random walk. 
\Qed

\begin{Definition}
The \defn{harmonic measure} for simple random walk on $\wdd$ starting at the top vertex $\e$ is the
probability measure on the real interval $[0,1)$ that is the law of the limit of the horizontal
\RL{positions} of the walk. 
Arguments similar to those used to prove Proposition \ref {p.hmexists} show that such a
measure, $\hmd$, exists.
\end{Definition}

Again, we may identify $\hmd$ explicitly as Haar (Lebesgue) measure on $\R/\Z$:

\begin{Proposition}\label{prop:dualharmonicuniform}
The harmonic measure $\hmd$ for the dual walk is uniform on $[0,1)$, equivalently, $\bigl(\frac{\delta_0 + \delta_1}2\bigr)^{\otimes \N}$.
\end{Proposition}
\begin{proof}
It suffices to show that for each positive integer $n$, this measure gives equal weight to the
intervals $\bigl[k2^{-n},(k+1)2^{-n}\bigr)$ for each $k$ between $0$ and $2^n-1$.

To see this equality, for each path $\rho$ converging to a point in one of these intervals, we pair
it with similar paths converging to points in each of the others. Define the dual leaving times
$\Td_n$ as the last time the walk on $\wdd$ is at depth $n$. For each $i$ between $0$ and $n-1$, at
time $\Td_i + 1$ the walk leaves depth $i$ for the last time. There are two \RL{equally likely} ways this could be done,
by appending $0$ or $1$. Consider the family of $2^n$ paths obtained by making each combination of
these $n$ \RL{independent} choices, and otherwise moving as in $\rho$. This grouping \RL{has the property}
that starting with any one of these $2^n$ paths produces the same set of $2^n$ paths. 

At each time $t > \Td_{n-1}$, these walks are at positions $\{x+\frac{l}{2^n}\}_{l=0}^{2^n-1}$ for some $x$. If one of these $2^n$ walks converges to a point $x$, then others converge to $x + \frac{l}{2^n}$ for each $l$. This completes the proof.
\end{proof}

\bigbreak
\section{Graph Symmetries and Probabilities}
\label{sec:reflections}
The graph $\w$ has only one nontrivial graph automorphism, described in Proposition \ref{prop:reflectionauto}. It interchanges the notions of `right' and `left' throughout the graph. 

\begin{Proposition}\label{prop:reflection}
If $\phi$ is the reflection automorphism of Proposition
\ref{prop:reflectionauto}, then the harmonic measure $\hm$ is
$\phi$-invariant.
\end{Proposition}
\begin{proof}
Automorphisms of any graph \RL{that fix the starting vertex} leave the law of simple random walk invariant.
\end{proof}

\begin{Corollary}\label{cor:onebitdensity}
If a number $x$ between $0$ and $1$ is chosen according to the harmonic measure $\hm$, then for any positive integer $n$, the probability that the $n$th binary digit of $x$ is a $0$ is $\frac{1}{2}$.
\end{Corollary}
\begin{proof}
As long as $x$ is not a dyadic rational number $\frac{k}{2^n}$, the map $\phi$ changes the $n$th bit of $x$. The probability that $x$ is such a dyadic rational is zero by Lemma \ref{lem:nopointmasses}, so Proposition \ref{prop:reflection} gives the result.
\end{proof}

Surprisingly, the generalisations of Corollary \ref{cor:onebitdensity} to strings of more than one bit are false.

\begin{Proposition}\label{prop:reflection2}
If a number $x$ between $0$ and $1$ is chosen according to the harmonic measure $\hm$, then for any positive integer $n$, the probabilities that the $n$th and $(n+1)$th binary digits of $x$ are $00$ or $11$ are equal, the probabilities of $01$ and $10$ are equal, and the former probabilities are \RL{strictly} greater than the latter.
\end{Proposition}

Intuitively, this is because \RL{if} the random walk starts in the right half of \rref f.wrapped-circ/, \RL{then} it is more likely to end up in the right half as well.
Before we prove Proposition \ref{prop:reflection2}, we give some preliminary results. The proof appears after Proposition \ref{prop:crestsum}.

First, we define a reflection on random walk paths, as in Proposition \ref{prop:reflection}. Unlike that reflection, this one will {\em not} be induced by an automorphism of the graph $\w$.

\begin{Definition}
Let $\rho$ be a path taken by the simple random walk on $\w$. Let $T$ be the first time greater than
the leaving time $T_0$ at which $\rho$ is at a string of the form $01 \conc 0^k$ or $11 \conc 0^k$
($k \ge 0$). It is possible that this never happens, in which case $T = \infty$.

Define a modified path $\phi_2(\rho)$ as follows. 
\begin{itemize}
\item If $T < \infty$, then $\phi_2(\rho)$ agrees with $\rho$ up until time $T$, and then proceeds as $\rho$ except with left and right moves switched.
\item If $T = \infty$, then $\phi_2(\rho) := \rho$.
\end{itemize}
\end{Definition}

We may regard $\phi_2$ as acting on the horizontal \RL{position} by a reflection in either 1/4 or 3/4, whichever is visited first; in fact, these reflections are
the same and implemented by the map $x \mapsto 1/2 - x \pmod 1$. The map $\phi_2$ cannot be derived from an automorphism of $\w$, because it sometimes interchanges the vertices $0$ and $1$, which have different degrees. However it only ever does this on a section of a path that never visits the vertex $\e$. Essentially, removing the vertex $\e$ increases the available set of automorphisms. %
The reader may wish to \RL{review} Figure \ref{fig:gamma0}.

\begin{Proposition}\label{prop:prer2pairing}
Let $\rho$ be a path that is not fixed by $\phi_2$ and that converges to $x$ in the interval $[0,1)$. Then $\phi_2(\rho)$ converges to $\frac{1}{2} - x$ modulo $1$. 
\end{Proposition}
\begin{proof}
If a sequence of horizontal \RL{positions} $x_i$ converges to $x$, then the sequence $(\frac{1}{2}-x_i)$ converges to
$\frac{1}{2} - x$.
\end{proof}

\begin{Corollary}\label{cor:r2pairing}
If a path $\rho$ is not fixed by $\phi_2$, then with probability $1$ exactly one of $\rho$ and $\phi_2(\rho)$ converges to a binary string starting with $00$ or $11$, while the other converges to a string starting with $01$ or $10$. 
\end{Corollary}
\begin{proof}
From Proposition \ref{prop:prer2pairing}, the map $\phi_2$ interchanges the open intervals $(0,\frac{1}{4})$ and $(\frac{1}{4},\frac{1}{2})$, in the sense that if $\rho$ converges to a point in $(0,\frac{1}{4})$ and $\phi_2(\rho) \neq \rho$, then $\phi_2(\rho)$ converges to a point in $(\frac{1}{4},\frac{1}{2})$, and vice versa. Likewise, $\phi_2$ interchanges the intervals $(\frac{1}{2},\frac{3}{4})$ and $(\frac{3}{4},1)$.

The statement of the present corollary differs from this result only in that it refers to the half-open intervals $[0,\frac{1}{4})$, $[\frac{3}{4},1)$, $[\frac{1}{4},\frac{1}{2})$, and $[\frac{1}{2},\frac{3}{4})$, and requires only that the claim be true with probability $1$. The probability that $\rho$ converges to any of the four points $0, \frac{1}{4}, \frac{1}{2}$, or $\frac{3}{4}$ is zero by Lemma \ref{lem:nopointmasses}, which completes the proof.
\end{proof}

We will prove Proposition \ref{prop:reflection2} by dividing paths into those fixed by $\phi_2$, and others. Corollary \ref{cor:r2pairing} shows that paths not
fixed by $\phi_2$ are as likely to converge to an element of the interval $(\frac{1}{4},\frac{3}{4})$ as to an element of the complement $(0,\frac{1}{4}) \cup (\frac{3}{4},1)$. It remains to consider paths that are fixed by $\phi_2$.

\begin{Proposition}\label{prop:r2fixing}
If $\rho$ is a path that is fixed by $\phi_2$ and that does not converge to $0, \frac{1}{4}, \frac{1}{2},$ or $\frac{3}{4}$, then there are two possibilities:
\begin{itemize}
\item At time $T_1$, $\,\rho$ is at $0$. Then after time $T_1$, $\,\rho$ only ever visits vertices beginning with $00$ or $11$, and so converges to a string starting with $00$ or $11$.
\item At time $T_1$, $\,\rho$ is at $1$. Then after time $T_1$, $\,\rho$ only ever visits vertices beginning with $01$ or $10$, and so converges to a string starting with $01$ or $10$.
\end{itemize}
\end{Proposition}
\begin{proof}
In the first case, at time $T_1 + 1$, the path $\rho$ is at $00$. It never returns to depth $1$, so the only way it could leave the set $[0,\frac{1}{4}) \cup (\frac{3}{4},1)$ is via sideways moves. However, this would result in it passing through $\frac{1}{4}$ or $\frac{3}{4}$. The path $\rho$ does not converge to $\frac14$ or $\frac34$, so it contains a sideways move after this time, which contradicts $\phi_2(\rho) = \rho$.

The second case is the same but with the vertex $00$ and the set $[0,\frac{1}{4}) \cup (\frac{3}{4},1)$ replaced by $10$ and $(\frac{1}{4},\frac{3}{4})$.
\end{proof}

We will also need to show that the set of walks that are fixed by $\phi_2$ has positive probability.

\begin{Proposition}\label{prop:maybefixed}
Consider a simple random walk on $\w$ that starts at $\e$. There is a positive probability that this walk eventually passes through the vertex $00$, never returns to depth $2$, and is fixed by $\phi_2$.
\end{Proposition}
\begin{proof}
It suffices to show that the probability of a walk starting at $00$ never reaching a vertex of the form $01 \conc 0^k$ or $11 \conc 0^k$ is positive.
This is immediate from Proposition \ref{p.hmexists}.
\end{proof}

These probabilities of last leaving depth 1 at either vertex 0 or 1 are related to hitting probabilities in a surprising way---we will see that they are equal to the probabilities that a random walk started at $0$ or $1$ ever reaches the vertex $\e$. Because the random walk leaves depth $1$ from either $0$ or $1$, these two probabilities sum to one. This will give the same relation for the hitting probabilities.

\begin{Definition}
Let $\pa(x \rightarrow y)$ be the set of paths from $x$ to $y$. If $x=y$, then this includes the path of length $0$. 
\end{Definition}

\begin{Definition}
If $\rho$ is a path, then $\Pr(\rho)$ is the probability of the path $\rho$---that is, the product of the probabilities of each step, which is equal to the product of the reciprocals of the degrees of each vertex, except for the final vertex. 
\end{Definition}

\begin{Definition}
The \defn{crest probability} from a vertex $\s$ of\/ $\w$ is the probability that the random walk,
started from $\s$, ever reaches the vertex $\e$. Denote this probability by $\esc(\s)$.
\end{Definition}

We now relate crest probabilities to the probabilities that a walk starting at $\e$ leaves depth $1$ for the last time at $0$ or at $1$, by reversing the paths in question.  

\begin{Remark}\label{rem:returnprob}
The probability $\esc(0)$ (respectively, $\esc(1)$) is equal to the probability that a walk starting at any vertex of degree $4$ (resp., $3$) ever reaches the level above its starting vertex.
\end{Remark}
\begin{proof}
Consider random walks starting at $0$ (resp., $1$) and at any other starting vertex of degree $4$ (resp., $3$), and couple them so that they always move in the same direction---up, left, right, or down. Either both will eventually reach the level above their starting vertices, or neither will. 
\end{proof}

\begin{Proposition} \label{prop:crestprobs}
The crest probabilities $\esc(0)$ and $\esc(1)$ are related by $$\esc(0) + \esc(1) = 1.$$
In fact, $\esc(i) = \P[X_{T_1} = i]$ for $i \in \{0, 1\}$ when $X_0 = \e$.
\end{Proposition}
\begin{proof}
By the craps principle \cite[p.~210]{Pitman}, $\P[X_{T_1} = i] = \P[X_{T_1} = i \mid T_0 = 0]$.
By Remark \ref{rem:returnprob}, $\P[X_{T_1} = i, T_0 = 0]$
is the sum of $\P(\rho) \cdot \frac1{\deg(i)} \cdot \bigl(1 - \esc(0)\bigr)$ over all paths $\rho$ that start at $\e$, end at $i$, and do not visit $\e$ again. Because $\P[T_0 = 0] = \frac1{\deg(\e)}\cdot\bigl(1 - \esc(0)\bigr)$, it follows that $\P[X_{T_1} = i]$ is the sum of $\P(\rho) \cdot \frac{\deg(\e)}{\deg(i)}$ over all paths $\rho$ that start at $\e$, end at $i$, and do not visit $\e$ again. This is the same as the sum of $\P(\rho)$ over all paths $\rho$ that start at $i$, end at $\e$, and do not visit $\e$ again, i.e., $\esc(i)$.
\end{proof}

These results generalise to the following.

\begin{Proposition}\label{prop:crestsum}
For any depth $n$, the sum of the crest probabilities from each of the $2^n$ vertices at depth $n$ is $1$. For each vertex $\s$ at depth $n$, the crest probability $\esc(\s)$ satisfies $\esc(\s) = \P[X_{T_n} = \s]$ when $X_0 = \e$.
\end{Proposition}
\begin{proof}
The proof is the same as that of Proposition \ref{prop:crestprobs}, with the vertex $i$ replaced by $\s$.
\end{proof}

\procl r.somerelns
Using symmetry and the condition that $\esc(\s)$ equals the average of the values of\/ $\esc$ at the neighbors of $\s$, one can show that $\esc(00) = 6\esc(0) - 3$, $\>\esc(10) = 3 - 5\esc(0)$, and $\esc(01) = \esc(11) = \esc(1)/2$.
\endprocl

We have a similar extension to walks on $\Gamma_+$:

\procl p.sideways
Let $\s$ be a vertex at depth $1$ of\/ $\Gamma_+$. The probability that simple random walk on $\Gamma_+$ started at $(0, \e)$ and conditioned never to visit depth $0$ last leaves depth $1$ at $\s$ equals the probability that simple random walk on $\Gamma_+$ started at $\s$ visits $(0, \e)$ before visiting any other vertex at depth $0$ (if any).
\qed
\endprocl

While these techniques relating crest probabilities to the positions at which a random walk last leaves an appropriate set could be applied to other graphs, we note that Remark \ref{rem:returnprob} requires that the graph in question be quite self-similar. 

\begin{proof}[Proof of Proposition \ref{prop:reflection2}]
Note that the probabilities that the first two bits of $x$ are $00$ or $11$ are equal, as are the probabilities to be $01$ or $10$, by the same argument as in Corollary \ref{cor:onebitdensity}. Thus it suffices to show that the probability of $(00 \text{ or } 11)$ is greater than that of $(01 \text{ or } 10)$.

Combining Corollary \ref{cor:r2pairing} with Proposition \ref{prop:r2fixing}, it suffices to show that the first case in Proposition \ref{prop:r2fixing} is strictly more likely than the second---that is, that among paths fixed by $\phi_2$, more of them leave depth $1$ for the last time from the vertex $0$ than from the vertex $1$.
Now the probability of a path being fixed by $\phi_2$ is independent of its location at time $T_1$. Therefore, the question reduces to showing that it is more likely to leave depth 1 from 0 than from 1. But these are the crest probabilities, for which this comparison is obvious because to return to $\e$ from 1 it is necessary to pass through 0, so $\esc(1) = \esc(0) \cdot \P_1\bigl[\pa(1 \rightarrow 0)\bigr]$, and $\P_1\bigl[\pa(1 \rightarrow 0)\bigr] \;<\; \P_1[T_1 > 0] \;<\; 1$, where $\P_1$ is the probability measure for the random walk starting at $1$. 
\end{proof}

\begin{Corollary} \label{c.hmnotequal}
The harmonic measure $\hm$ and the dual harmonic measure $\hmd$ are not equal. 
\end{Corollary}
\begin{proof}
Proposition \ref{prop:reflection2} shows that $\hm$ is not uniform on $[0,1)$, and Proposition \ref{prop:dualharmonicuniform} says that $\hmd$ is uniform on $[0,1)$.
\end{proof}

Not only are these two measures not equal to one another, they are mutually singular.

\begin{Proposition}\label{prop:mutuallysingular}
The harmonic measure $\hm$ and dual harmonic measure $\hmd$ are mutually singular.
\end{Proposition}
\begin{proof}
The harmonic measure $\hm$ is invariant under the left shift; in fact, we need a specific form of $\hm$ showing this invariance. Write $(b_j)_{j \ge 1} \in \{0, 1\}^{\Z^+}$ for the limit of the random walk
$(X_t)_{t \ge 0}$ on $\w$. Write $(M_t)_{t \ge 0}$ for the sequence of steps from $\Sigma :=
\{\mathrm L, \mathrm R, \mathrm U, \mathrm D\}$ taken by the random walk.
Thus, $X_{t+1}$ is obtained from $X_t$ by applying the move $M_{t+1}$.
Then there is a measurable function $f \colon \Sigma^{\Z^+} \to \{0, 1\}$ such that $b_1 =
f\bigl((M_{T_0+t})_{t \ge 1}\bigr)$.
In fact, the same function $f$ gives all bits as $b_j = f\bigl((M_{T_{j-1}+t})_{t \ge 1}\bigr)$.
The sequences $(M_{T_j + t})_{1 \le t \le T_{j+1} - T_j}$ are IID for $j \ge 0$, as noted in the
proof of \rref t.stationary/. 
Therefore, $(b_j)_{j \ge 1}$ is a factor of this IID sequence, so its law, $\hm$,
is ergodic, as is, obviously, $\hmd$. Two ergodic measures for the same
transformation are either equal or mutually singular, whence the result.
\end{proof}

We remark that because $\hm$ is a factor of IID, it is isomorphic to what is called in ergodic theory
a Bernoulli shift.

The stationary measure, by contrast, is not even shift-invariant: 
\eqaln
{\sm\bigl\{(a_k)_{k \le 0} : a_0 = 1\bigr\} 
\;=\;
p_3 \;\le\; 3/7 &\;<\; 1/2 
\;=\; \hm\bigl([1/2, 1)\bigr) 
\\ &\;=\; \lim_{j \to - \infty} \sm\bigl\{(a_k)_{k \le
0} : a_j = 1\bigr\}}
by \rref p.smhm/.

We may now answer Curien's Question \ref{q.curien}. Recall that $\rstates$ and $\rstatesd$ have a common coding by $\Z_2$.

\procl t.curien
There is no stationary measure on $\rstates$ that induces a measure on $\rstatesd$ that is absolutely
continuous with respect to some stationary measure. 
\endprocl

\rproof
Both stationary measures are unique, so it suffices to show that the stationary measure $\sm$ on $\rstates$ and the stationary measure on $\rstatesd$ are mutually singular. Proposition \ref{Proposition:smhm} shows that $\sm$-almost all binary strings have substring densities given by $\hm$; the same is evident for $\smd$ and $\hmd$ by Propositions \ref{Proposition:sdmuniform} and \ref{prop:dualharmonicuniform}.  Corollary \ref{c.hmnotequal} says that $\hm$ and $\hmd$ are not equal, so $\sm$ and $\smd$ are mutually singular.
\Qed

\bigbreak
\section{Numerical Results} \label {sec:numerics}

Because we do not have exact expressions for either $\sm$ or $\hm$, we devote this section to numerical approximations. Some interesting patterns will emerge, leading to some open questions.

\subsection{Numerical bounds on \texorpdfstring{$p_3$}{p\textunderscore 3}}

One of the most interesting quantities is
$p_3$, the proportion of time spent at vertices of degree $3$, well defined by Proposition \ref{prop:p3exists}. We do not know the exact value of $p_3$, but the following sections give bounds on this quantity. Rough analytic bounds were obtained in Remark \ref{r.moreonspeed}.

We begin with a discussion of the entire stationary measure, $\sm$.
An approximation of $\sm$ is shown in \rref f.statmsr/. In order to show
$\sm$, we map $\Z_2 \to [0, 1]$ by $\sum_{k \ge 0} a_k 2^k \mapsto \sum_{k \ge 0} a_k 2^{-k-1}$ for
$a_k \in \{0, 1\}$. We then push forward $\sm$ via this map. We approximated $\sm$ by using the
corresponding Markov chain on binary strings of length 12; removing the first bit corresponding to a
move upwards entailed adding a last bit, which we fixed to be 0. We aggregated the stationary measure
for this finite system into blocks of the first 6 most significant bits. Multiplying each of those
numbers by $2^6$ shows it as a density in \rref f.statmsr/.
While we do not know how accurate these probabilities are, they appear surprisingly accurate if we can judge by the implied aggregate estimate of $p_3$, which would be about $0.382332$; this differs by only $10^{-6}$ from our best estimate of $p_3$ using another method explained in the next two paragraphs. In any case, one can show that as the length used for the binary strings in this finite system grows, the error decays exponentially in that length.
Also, because we used always 0 when a new last bit was needed, one might expect that the resulting estimate of $p_3$ would be a lower bound; however, we do not know a proof of this.

\begin{figure}[ht]
\begin{center}
\includegraphics[width=.7\textwidth]{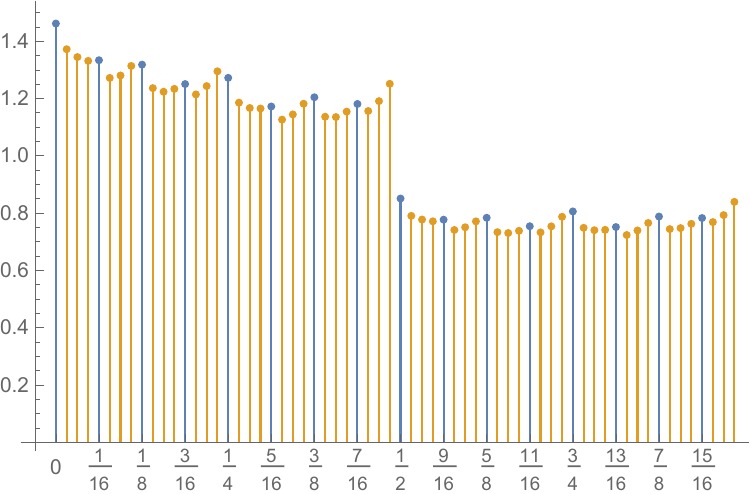}
\end{center}
\caption{A finite-size approximation of the stationary measure $\sm$ for the distribution of the
first 6 bits, given as a density.}
\label{f.statmsr}
\end{figure}

Notice that the speed $p_3/3$ is the probability of never returning to the current level, because each level is left for the last time just once. On the other hand, this nonreturn probability is equal to $\bigl(p_3/3+(1-p_3)/4\bigr)\bigl(1-\esc(0)\bigr)$. Equating these two expressions for the speed, we obtain
\[
\esc(0) = \frac{3(1-p_3)}{3+p_3}
\quad \text{ and }
\quad
p_3 = \frac{3\bigl(1-\esc(0)\bigr)}{3+\esc(0)}.
\]
Using this and \RL{an} estimate $\esc(0) \approx 0.547846 \pm 10^{-6}$, our estimate for $p_3$ is
$0.382333 \pm 10^{-6}$.

\RL{Our} estimate for $\esc(0)$ was obtained by the following method.
The function $\s \mapsto \esc(\s)$, defined on the set of vertices $\s$ of $\w$, is the solution to the Dirichlet problem with $\esc(\e) = 1$ and $\lim_{\depth(\s) \to\infty} \esc(\s) = 0$. In other words, if $\w(n)$ is the graph induced by vertices in $\w$ whose depth is between 0 and $n$ and $\esc_n$ is the harmonic function on $\w(n)$ with $\esc_n(\e) = 1$ and $\esc_n(\s) = 0$ for all $\s$ of depth $n$, then $\lim_{n \to\infty} \esc_n = \esc$ pointwise. Here, `harmonic' means that $\esc_n(\s)$ is the average value of $\esc_n$ at the neighbours of $\s$ for every $\s$ whose depth is between 1 and $n-1$. In addition, $\esc_n$ increases in $n$. The function $\esc_n$ solves a sparse linear system of equations, which we solved for $2 \le n \le 20$.
We show $\log_2 \esc_{20}(\s)$ for $\s$ of depth at most 12 in \rref f.crestprs/. There is a clear pattern based on the least significant bits; this is explained by Remark \ref{rem:returnprob}.
Despite Proposition \ref{prop:crestsum}, the sum $\esc_n(0) + \esc_n(1)$ is not 1, but only tends to 1 as $n \to\infty$. Thus, we normalised $\esc_n$ to sum to 1 on each level. The last 14 of the resulting numbers $\esc_n(0)$ seemed to be approaching a limit exponentially fast with ratio 2, so we fit such a curve to them, leading to our estimate of the preceding paragraph. One can show that $\esc_n(0)$ does approach $\esc(0)$ exponentially fast; one can also get an upper bound on $\esc(0)$ by solving the Dirichlet problem where the values $\esc_n(\s)$ for $\s$ of depth $n$ are set to an upper bound on $\esc(\s)$ for such $\s$. 

\begin{figure}[ht]
\begin{center}
\includegraphics[width=\textwidth]{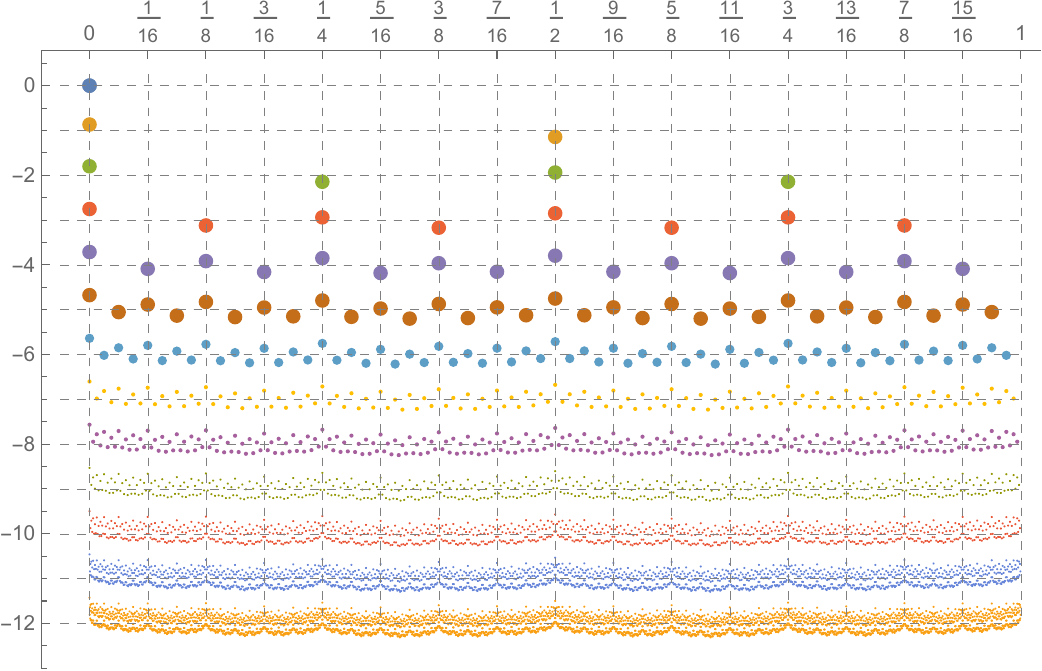}
\end{center}
\caption{Approximation of the base-2 logarithms of the crest probabilities at depths 0--12. The dot sizes vary only for visibility.}
\label{f.crestprs}
\end{figure}

\subsection{Exploration of the harmonic measure}

In this section, we investigate further the harmonic measure, $\hm$, for simple random walk on $\w$. We know from Proposition
\ref{prop:mutuallysingular} that $\hm$ is singular with respect to Lebesgue measure. Because $\hm$ is
the quotient of $\hp$ by $\Z$, it follows that $\hp$ is also singular with respect to Lebesgue
measure on $\R$.

Figure \ref{fig:whm} shows an approximation to $\hm$.
The figure has one dot for each interval of length $2^{-14}$, whose ordinate denotes the measure of that interval times $2^{14}$, as if it were a density.
Note that because $\hm$ is singular, finer approximations would tend to 0 and $\infty$ a.e.\ with respect to Lebesgue measure.
We calculated this approximation as follows. Consider simple random walk on $\Gamma_+$ starting at depth 0. Let $K_n$ be the change in horizontal position between times $T_{n-1}$ and $T_n$ for $n \ge 1$. The reasoning behind Proposition \ref{prop:drift} shows that $(2^n K_n)_{n \ge 1}$ are IID. Note that $2^n K_n \in \Z$. The distribution of $\sum_{n \ge 1} K_n$ is $\hp$; taking this modulo 1 yields $\hm$. Thus, it suffices to know the law of $K_1$. By \rref p.sideways/, $\P[K_1 = k + i/2]$ equals the probability that simple random walk starting at $(k, i)$ at depth $1$ visits $(0, \e)$ before visiting any other vertex at depth $0$, where $k \in \Z$ and $i \in \{0, 1\}$.
This is a solution to a Dirichlet problem again; we approximated it by the corresponding Dirichlet problem on $\w$ between depths $6$ and $19$, finding the probability for each $\s$ of depth $7$ that random walk from $\s$ visits $0^6$ before visiting any other vertex at depth $6$ or any vertex at depth $19$. We normalised these probabilities to add to 1. Having this approximation to the law of $K_1$ at hand, we approximated the law of $\sum_{n = 1}^{22} K_n \pmod 1$ and then
aggregated to intervals of length $2^{-14}$.

\begin{figure}[ht]
\begin{center}
\includegraphics[width=.9\textwidth]{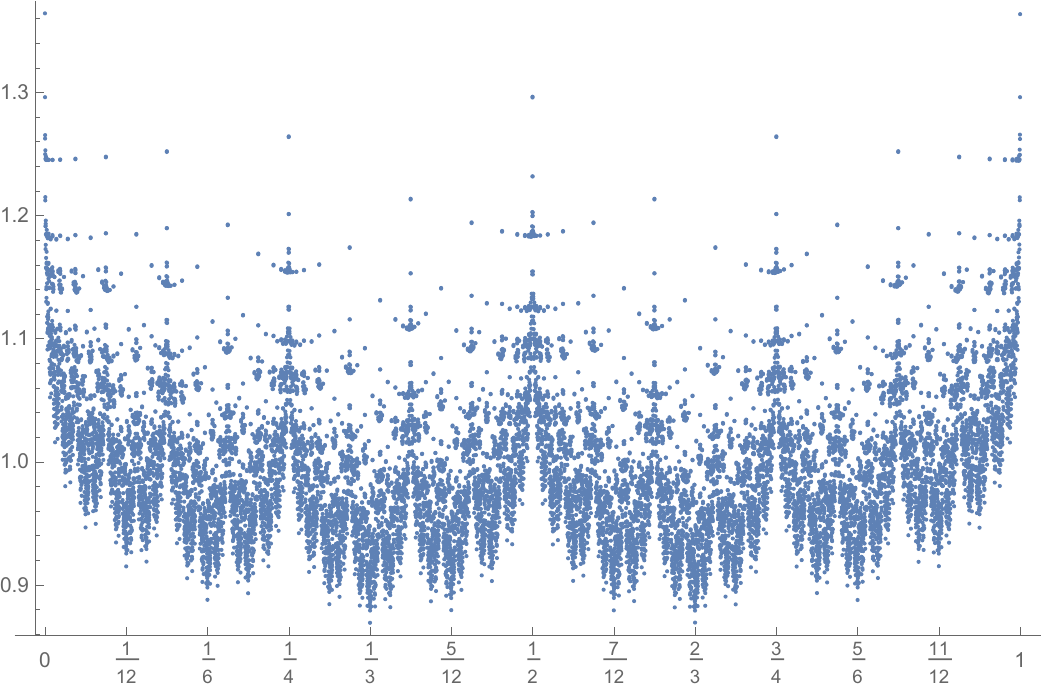}
\end{center}
\caption{An approximation of the harmonic measure $\hm$ on the interval $[0,1]$, using intervals of
length $2^{-14}$.} 
\label{fig:whm}
\end{figure}

\begin{figure}[ht]
\begin{center}
\includegraphics[width=.9\textwidth]{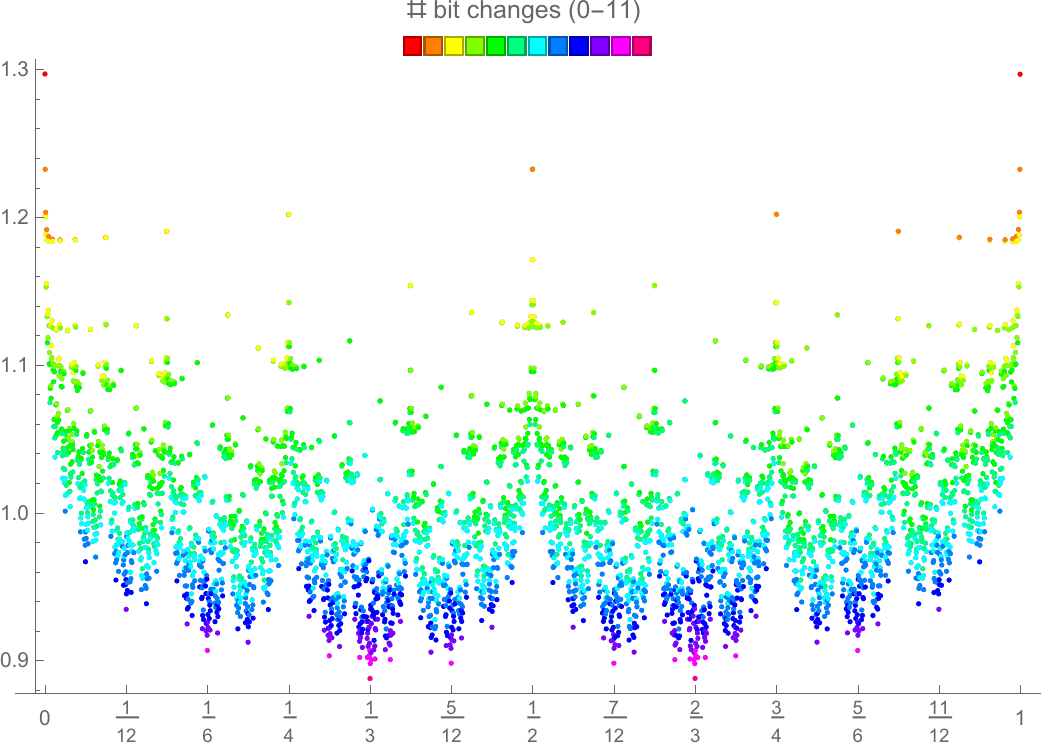}
\end{center}
\caption{An approximation of the harmonic measure $\hm$ on the interval $[0,1]$, using intervals of
length $2^{-12}$. The colour shows the number of differences in successive bits.} 
\label{f.bitchanges}
\end{figure}

The locations of the most extreme maxima and minima in Figure \ref{fig:whm} appear to be controlled by binary representations. Maxima occur at positions whose binary expressions are short and terminate, like $0, \frac{1}{2}, \frac{1}{4}, $ and $\frac{3}{4}$. Minima occur at positions whose binary expressions end in alternating sequences of $0$'s and $1$'s, like $\frac{1}{3}$, $\frac{2}{3}$, and $\frac{1}{6}$.
\rref f.bitchanges/ shows \RL{the same plot as} Figure \ref{fig:whm}, but with intervals of length $2^{-12}$ instead of $2^{-14}$ \RL{and}
with colours corresponding to the number of
differences in successive bits in the string of length 12.
One can explain such a relationship by the use of $g$-measures; see below.

While $\hm$ is self-similar in that it is invariant under the map $x \mapsto 2x \pmod 1$,
Figure \ref{fig:whm} appears to show a different kind of self-similarity---the portion of the graph between $0$ and $0.5$ exhibits a similar pattern of oscillations to the whole graph.
This indicates a weak influence of the first bit on the distribution of the remaining bits.
We look at this quantitatively next.

For $x \in [0, 1)$ and $i \in \{0, 1\}$, consider the conditional probability  $p(i, x) := \hm\bigl\{\lfloor{2x}\rfloor = i \mid 2x
\pmod 1\bigr\}$.
This is defined for $\hm$-a.e.\ $x$. If we write $x$ as a binary string $(x_k)_{k \ge 1}$, then $p(i,
x)$ is the $\hm$-probability that $x_1 = i$ given $(x_k)_{k \ge 2}$.
The base-2 entropy of $\hm$ equals 
\[
h := \int_0^1 \bigl[- p(0, x) \log_2 p(0, x) - p(1, x) \log_2 p(1, x)\bigr] \,d\hm(x)
=
\int_0^1 - \log_2 g(x) \,d\hm(x),
\]
where $g(x) := p\bigl(\lfloor{2x}\rfloor, x\bigr)$ is defined $\hm$-a.e.
The entropy $h$ is also the Hausdorff dimension of $\hm$: there is a set of dimension $h$ that carries $\hm$,
but no set of dimension smaller than $h$ carries $\hm$ \cite {Billingsley:info}.
\RL{Because $\hm$ is invariant under the map $x \mapsto 2x \pmod 1$, yet $\hm$ is not
equal to Lebesgue measure, it follows that $h < 1$.}
Figure \ref{fig:whmratio} shows an approximation to the function $g$
using $2^{12}$ points. This calculation used the approximation to $\hm$ mentioned above.
Using our approximations of $g$ and $\hm$, this gives $h \approx 0.999799$. However, entropy is notoriously difficult to estimate, so although we have other methods that support this estimate \RL{of $h$}, we cannot claim great confidence in it.
Interestingly, Figure \ref{fig:whmratio} appears to show a continuous curve,
monotone on each half of the interval. 
If such a
curve really does exist, then harmonic measure is called a $g$-measure
\cite{Keane}.
This curve appears to have bounded variation with its derivative being a singular measure: see Figure \ref{fig:whmratio-deriv}.
Proposition \ref{prop:reflection} shows that
$g(x) = g(1-x)$. Also, $g(x) + g(x+1/2) =
1$ by definition.
We believe that not only is $g$ continuous and monotone decreasing on $[0, 1/2]$, but that it determines $\hm$ uniquely as a $g$-measure; see \cite{Ledrappier,Stenflo,BramsonKalikow} for
discussions of uniqueness.

Assuming that $g$ is continuous and appears as in Figure \ref{fig:whmratio}, we have that $g(x) > 1/2$ for $x \notin [1/4, 3/4]$ and $g(x) < 1/2$ for $x \in (1/4,
3/4)$. Therefore, $\hm$ will tend to be large at $x$ with $2^k x \notin [1/4, 3/4] \pmod 1$ for most $k$ (the same as the first two bits of $2^k x$ being the
same), i.e., for $x$ with few changes in successive bits, and $\hm$ will tend to be small at $x$ with $2^k x \in (1/4, 3/4) \pmod 1$ for most $k$ (the same as the
first two bits of $2^k x$ being different), i.e., for $x$ with many changes in successive bits. This would explain \rref f.bitchanges/.

\begin{figure}[ht]
\begin{center}
\includegraphics[width=0.9\textwidth]{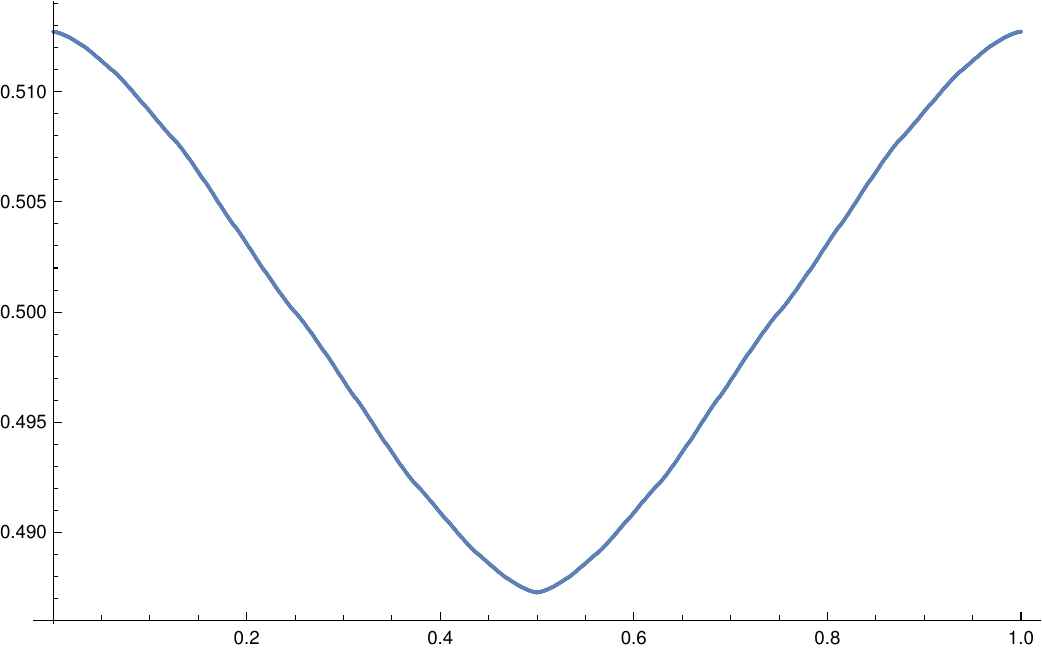}
\end{center}
\caption{The $\hm$-probability $g(x)$ of the first bit of $x$ given the rest of the bits of $x$, apparently continuous in $x$ and close to 1/2 and taking the value
exactly 1/2 at $x \in \{1/4, 3/4\}$. The plot shows $g(x)$ for $x$ a multiple of $2^{-12}$.} 
\label{fig:whmratio}
\end{figure}

\begin{figure}[ht]
\begin{center}
\includegraphics[width=0.45\textwidth]{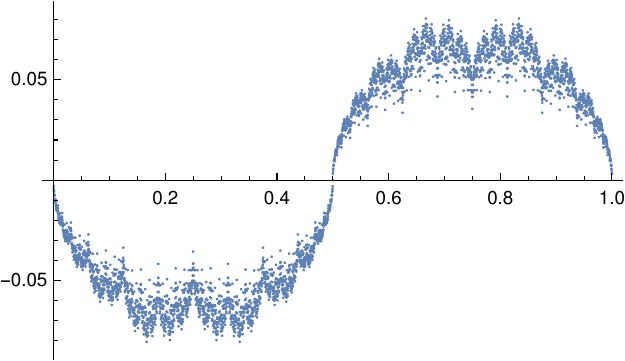}
\hskip20pt
\includegraphics[width=0.45\textwidth]{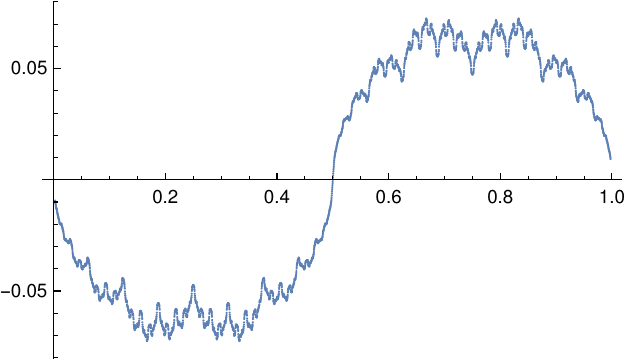}
\end{center}
\caption{The discrete derivative of the curve in Figure \ref{fig:whmratio} with step $2^{-12}$ on the left and with step $21 \cdot 2^{-12}$ on the right. This latter quantity is chosen to illustrate that the picture is similar for larger step sizes, and it appears similar for all step sizes in this range, only gradually losing finer details.}
\label{fig:whmratio-deriv}
\end{figure}

\bigbreak
\section{Optimality of the Construction}\label{sec:optimality}

In this section, we show that any graph which is a counterexample to our main question must have either some vertices of degree at least $5$ or some faces of degree at
least $5$. Our graphs $\Gamma_a$ have vertices of degrees $3$ and $4$ and faces of degree $5$.

\begin{Proposition}\label{prop:fewdegree3}
If an infinite plane graph has vertices only of degree at most $4$ and faces only of degree at most $4$, then the number of vertices of degree $3$ is at most $4$.
\end{Proposition}
\begin{proof}
We use the combinatorial curvature of a vertex, defined to be $1$ minus half of the degree plus the sum of the reciprocals of the degrees of the incident
faces. By Corollary 1.4 of \cite{DevosMohar}, the sum of the combinatorial curvatures over all vertices is at most $2$, and at most $1$ if the graph is
infinite. 

Write the combinatorial curvature at a vertex $v$ as $1$ plus the sum over all faces $f$ incident to $v$ of $1/\deg f - 1/2$.
The assumptions that vertex and face degrees are at most $4$ imply that the curvature at each vertex is at least $0$, with equality if and only if the vertex has degree $4$ and all its incident faces are squares. The curvature is at least $\frac{1}{4}$ for a vertex of degree $3$.  This shows that there are at most $4$ vertices of degree $3$.
\end{proof}

\begin{Corollary}\label{cor:fewtriangles}
If an infinite plane graph has vertices only of degree at most $4$ and faces only of degree at most $4$, then the number of vertices of degree $3$ plus the number of triangular faces is at most $4$.
\end{Corollary}
\begin{proof}
Apply Proposition \ref{prop:fewdegree3} to the graph formed from the primal and dual together, with new vertices of degree $4$ where the primal and dual edges cross. This new graph has one vertex of degree $3$ for each primal vertex of degree $3$ and one vertex of degree $3$ for each triangular face.
\end{proof}

The upper bound of $4$ is sharp, as can be seen by taking a triangular cylinder formed of squares that is infinite in one direction and capped by a triangle at the other.

\begin{Corollary}
The only stationary infinite plane graph whose vertices and faces have degree at most $4$ is the square lattice.
\end{Corollary}
\begin{proof}
Corollary \ref{cor:fewtriangles} shows that such a graph can have only finitely many vertices and faces whose degree is not $4$. If it is to be stationary, then it must have either zero or infinitely many such vertices and faces, so all vertices and faces have degree $4$. In this case, the graph is the square lattice.
\end{proof}

\bibliographystyle{plain} 
\bibliography{bib} 

\end{document}